%% file: TCNS-19-0006.final.arXiv.tex
\documentclass[final]{IEEEtran}
\usepackage{amsmath,amssymb,amsthm,epsfig,dsfont,color,subfigure,empheq,graphicx}
\usepackage{enumerate,url,algpseudocode,algorithm,wasysym,epstopdf,verbatim,array,textcomp,balance}
\usepackage[table]{xcolor}

\input{stylesV2}

\begin{document}
\title{Optimal Scheduling of Water Distribution Systems}

\author{Manish K. Singh~\IEEEmembership{Student Member,~IEEE} and Vassilis Kekatos,~\IEEEmembership{Senior Member,~IEEE}

\thanks{M. K. Singh and V. Kekatos are with the Bradley Dept. of ECE, Virginia Tech, Blacksburg, VA 24061, USA. Emails:\{manishks,kekatos\}@vt.edu. Work partially supported by the US National Science Foundation grant 1711587.}
}

\markboth{IEEE TRANSACTIONS ON CONTROL OF NETWORK SYSTEMS (to appear)}{Singh and Kekatos: Optimal Scheduling of Water Distribution Systems}

\maketitle
\begin{abstract}
With dynamic electricity pricing, the operation of water distribution systems (WDS) is expected to become more variable. The pumps moving water from reservoirs to tanks and consumers, can serve as energy storage alternatives if properly operated. Nevertheless, optimal WDS scheduling is challenged by the hydraulic law, according to which the pressure along a pipe drops proportionally to its squared water flow. The optimal water flow (OWF) task is formulated here as a mixed-integer non-convex problem incorporating flow and pressure constraints, critical for the operation of fixed-speed pumps, tanks, reservoirs, and pipes. The hydraulic constraints of the OWF problem are subsequently relaxed to second-order cone constraints. To restore feasibility of the original non-convex constraints, a penalty term is appended to the objective of the relaxed OWF. The modified problem can be solved as a mixed-integer second-order cone program, which is analytically shown to yield WDS-feasible minimizers under certain sufficient conditions. Under these conditions, by suitably weighting the penalty term, the minimizers of the relaxed problem can attain arbitrarily small optimality gaps, thus providing OWF solutions. Numerical tests using real-world demands and prices on benchmark WDS demonstrate the relaxation to be exact even for setups where the sufficient conditions are not met. 
\end{abstract}

\begin{IEEEkeywords}
Water flow equations, convex relaxation, second-order cone constraints, optimal water flow.
\end{IEEEkeywords}

\section{Introduction}
While WDS serve as a critical infrastructure, there is an increasing emphasis on improving their reliability, quality, and efficiency. The cost-intensive installation and maintenance of WDS components, such as pipelines, pump stations, and reservoirs, have motivated network planning studies~\cite{Sherali1998enhanced}, \cite{Sherali2001}, \cite{zong2006harmony}, \cite{bragalli2015mathprog}. From an operational perspective, a recent survey on WDS optimization identifies pump scheduling and water quality as the two focus areas~\cite{Helena2017lost}. Recognizing that 4\% of the total electricity consumption in the United States is attributed to water network operations~\cite{David2008fourpercent}, and that the electricity cost for pumping constitutes the largest expenditure for water utilities~\cite{Vanzyl2004}, stresses the significance of optimal WDS scheduling. 

A typical WDS schedule would run pumps mainly at night when electricity prices are low to transfer water from reservoirs through pipes and fill up elevated tanks located closer to water demands. Under the smart city vision, dynamic electricity pricing and demand-response programs incentivize more flexible WDS schedules to minimize operational costs. For example, a surplus of residential solar generation around midday could be locally consumed to run pumps and fill up pumps, thus serving as an energy storage alternative. Adaptive WDS scheduling and the anticipated joint dispatching of electric power and water networks, motivate the need for scalable optimization tools and more realistic system models. 

The operation of WDS is constrained by minimum pressure requirements; capacity limitations imposed by pumps, pipelines, and tanks; and a set of hydraulic constraints. It is exactly these hydraulic constraints that give rise to complex mixed-integer and nonlinear formulations, and have been dealt so far in three broad ways~\cite{Helena2017lost}. The first class of methods enforces pressure and capacity constraints explicitly, while the hydraulic constraints are included implicitly through water network simulation tools, such as EPANET \cite{epanet2000}, \cite{odan2015epanet}. Metaheuristic approaches such as genetic algorithms \cite{Vanzyl2004}, ant-colony optimization \cite{hashemi2014ant}, or limited discrepancy search \cite{Ghaddar2015lagrange}, are then used together along with a WDS simulator to obtain an operating point. Some variants replace the slow but exact simulator with surrogate WDS models based on artificial neural networks or interpretive structural models \cite{broad2010meta}, \cite{Arai2013ism}. It has been demonstrated however that WDS optimization using metaheuristics coupled with a simulator scales unfavorably due to the computational effort required~\cite{Giacomello2013hybrid}.

The second class of methods rely on formulating (mixed-integer) nonlinear programs and handling them via nonlinear solvers~\cite{IBM2012valve}. 
A mixed-integer second-order cone formulation for optimal pump scheduling relaxes the hydraulic constraints to render the problem convex in the continuous variables~\cite{TaylorOWFCDC}, \cite{TaylorOWFcones}. The relaxation is shown to be exact presuming all pipes are equipped with pressure-relieving valves and upon ignoring some pressure tank constraints. The water-power nexus has been studied in \cite{Parvania2017GM}, wherein the non-convex hydraulic constraints are passed on to a non-convex solver with no optimality guarantees. The security of interdependent water-power-gas networks has been studied from a game-theoretic viewpoint in \cite{Saad2017game}, using the non-convex hydraulic constraints. 

The third class of methods uses linearization to end up with a computationally tractable mixed-integer linear program (MILP) formulation~\cite{bragalli2015mathprog}, \cite{Verleye_tankpumpmodel}. Adopting~\cite{TaylorOWFcones} to find an optimal water-power flow dispatch, reference~\cite{ZamzamOWPF} handles the non-convex constraints arising from both water and electric power networks via a successive convex approximation technique. The latter approach features computational advantages without the inaccuracies of linearization; yet water flow directions and the on/off status of pumps are assumed given. The participation of WDS in demand response and frequency regulation through pump scheduling with piece-wise linearization of hydraulic constraints has been suggested in \cite{Parvania2018DR}, \cite{Parvania2018WatPower}, \cite{Menke2017DR}. 

Towards computationally convenient WDS solvers, the contribution of this work is two-fold. First, a generalized model for various WDS components is developed in Section~\ref{sec:model}. Some of its distinct features include separability of binary and continuous variables, flexibility of bypassing pumps, bidirectional flows, and precise modeling of tank operation. Second, an OWF problem to minimize electricity operation cost for fixed-speed pumps is put forth in Section~\ref{sec:PF}. Sections~\ref{sec:CR}--\ref{sec:PCR} develop a convex relaxation, which is later augmented by a novel penalty term to promote minimizers that are feasible for the water network. Under specific conditions, the penalized relaxation is shown to yield a minimizer of the original non-convex OWF problem. The numerical tests of Section \ref{sec:tests} on benchmark WDS corroborate that the proposed relaxations can yield feasible and optimal WDS dispatches even when the analytical conditions are grossly violated.

\section{Water Network Modeling}\label{sec:model}
A water distribution system can be represented by a \emph{directed} graph $\mcG:=(\mcM,\mcP)$. Its nodes indexed by $m \in \mcM$ correspond to water reservoirs, tanks, and points of water demand. Reservoirs serve as primary water sources and constitute the subset $\mcM_r \subset \mcM$. Similarly, the nodes hosting tanks comprise the subset $\mcM_b\subset \mcM$. The nodes in $\mcM_r \cup \mcM_b$ do not serve water consumers. This is without loss of generality, since a potential co-located consumer at a node $m\in \mcM_r\cup \mcM_b$ can be attached to an auxiliary node connected to the node $m$ through a lossless pipe. Let $d_m^t$ be the rate of water injected into the WDS from node $m$ during period $t$. Apparently, for reservoirs $d_m^t\geq 0$; for demand nodes with water consumers $d_m^t\leq0$; tanks may be filling or emptying; and for junction nodes $d_m^t=0$. 

The elements of the edge set $\mcP$ of $\mcG$ represent water pipes, and their cardinality is $P:=|\mcP|$. All edges in $\mcP$ are assigned an arbitrary direction. The \emph{directed} edge $(m,n)\in\mcP$ models the pipeline linking nodes $m$ and $n$. If $(m,n)\in\mcP$, then $(n,m)\notin\mcP$. The water flow on edge $(m,n)$ is denoted by $d_{mn}^t$. If water runs from node $m$ to node $n$ at time $t$, then $d_{mn}^t\geq0$; and negative, otherwise. Flow conservation dictates
\begin{equation}\label{eq:node1}
d_m^t=\sum_{k:(m,k)\in\mcP} d_{mk}^t - \sum_{k:(k,m)\in\mcP} d_{km}^t,\quad \forall m,t.
\end{equation}

In addition to water injections and flows, water distribution system (WDS) operation is also governed by pressures. Water pressure is typically surrogated by the quantity of pressure head, which is measured in meters and is linearly related to water pressure \cite{Verleye_tankpumpmodel}. In detail, a pressure head of $h$ meters corresponds to a water pressure of $h\rho \tilde{g}$ pascal, where $\rho$ is the water density in kg/m$^3$, assumed to be a known constant and $\tilde{g}$ is the acceleration due to gravity in m/sec$^2$. The pressure head (also known as piezometric pressure head) at a node equals its geographical elevation plus the manometric pressure head attributed to the height of the water column or pumps.

\begin{table}[t]
\renewcommand{\arraystretch}{1.3}
\caption{Nomenclature}
\vspace*{-1em}
\begin{center}
\label{tbl:nomenclature}
\begin{tabular}{|l|l|}
\hline \hline 
\textbf{Symbol} & \textbf{Meaning}\\
\hline 
\hline $\mcM$  & node set\\
\hline $\mcM_r,\mcM_b$ & node sets of reservoirs and tanks\\
\hline $\mcP$, $P$ & edge set and number of edges\\
\hline $\mcP_a,\bmcP_a$ & edge set hosting pumps and its complement\\
\hline $d_m^t$ & injection at node $m$ and time $t$\\
\hline $d_{mn}^t$& flow on edge $(m,n)$ at time $t$\\
\hline $\tilde{d}_{mn}^t,\underline{d}_{mn},\overline{d}_{mn}$& flow through pump $(m,n)$ at time $t$, and limits\\
\hline $\bd$ $(\tbd)$ & pipe (pump) flows at all times\\
\hline $h_m^t,\underline{h}_m$& pressure at node $m$ during time $t$, and limit\\
\hline $\bh$ & nodal pressures at all times\\
\hline $c_{mn}$ & loss (consumption) coefficient of pipe (pump)\\
\hline $x_{mn}$ & flow direction (running status) for pipe (pump)\\
\hline $g_{mn}$ & pressure added by pump\\
\hline $\alpha_m^t$ & connectivity status for reservoir or tank $m$\\
\hline $\bar{h}_m$ & constant pressure at reservoir $m$\\
\hline $\beta_m^t$ & filling/emptying status of tank $m$\\
\hline $\ell_m^t,\underline{\ell}_m,\overline{\ell}_m$ & water level in tank at time $t$ and its limits\\
\hline $A_m$ & cross-sectional area for tank $m$\\
\hline $\delta$ & time interval\\
\hline $\pi_t$ & electricity cost at time $t$\\
\hline $f(\tbd)$ & total pumping cost given pump flows $\tbd$\\
\hline $\bA(\bd^t)$ & incidence matrix based on flow directions at $t$\\
\hline $g(\bh)$ & penalty function\\
\hline $\lambda$ & weighing parameter\\
\hline $M$ & big-$M$ trick parameter\\
\hline \hline
\end{tabular}
\end{center}
\end{table}

The pressure head or henceforth simply \emph{pressure} at node $m$ during time $t$ will be denoted by $h_m^t$. The operation of water networks requires a minimum manometric pressure at all nodes $m$. Adding this common minimum value of manometric pressure to the specific but known geographical elevation of each node $m\in\mcM$ gives a lower limit on its pressure as
\begin{equation}\label{eq:hmin}
h_m^t\geq \underline{h}_m.
\end{equation} 
Water movement in a pipe results in a quadratic pressure drop. In detail, the pressure drop across pipeline $(m,n)\in \mcP$ is described by the Darcy-Weisbach equation \cite{Verleye_tankpumpmodel}
\begin{equation}\label{eq:headloss}
h_m^t-h_n^t=c_{mn}\sign(d_{mn}^t)(d_{mn}^t)^2
\end{equation}
where the loss coefficient $c_{mn}:=\frac{\ell_{mn}f_{mn}}{4\pi^2r_{mn}^5\tilde{g}}$ depends on the pipe length $\ell_{mn}$; its inner radius $r_{mn}$; and the Darcy friction factor $f_{mn}$. Although factor $f_{mn}$ actually depends on flow $d_{mn}$ in a continuous nonlinear manner, it is typically approximated as constant; see~\cite{bragalli2015mathprog} and references therein. The $\sign$ function is defined such as $\sign(0)=0$ and it ensures that pressure drops in the direction of water flow. To avoid the discontinuity of the sign, we propose a mixed-integer model using the big-$M$ trick for the pressure drop in pipeline $(m,n)$ using the binary variables $\{x_{mn}^t\}_{t=1}^T$. In particular, the pressure drop equation of \eqref{eq:headloss} can be equivalently expressed through the constraints
\begin{subequations}\label{eq:headlossMI}
\begin{align}
-M(1-x_{mn}^t)&\leq d_{mn}^t\leq Mx_{mn}^t\label{seq:hla}\\
-M(1-x_{mn}^t)&\leq h_m^t-h_n^t-c_{mn}(d_{mn}^t)^2\leq M(1-x_{mn}^t)\label{seq:hlb}\\
-Mx_{mn}^t&\leq h_m^t-h_n^t + c_{mn}(d_{mn}^t)^2\leq Mx_{mn}^t\label{seq:hlc}\\
x_{mn}^t&\in\{0,1\}
\end{align}
\end{subequations}
for a large $M>0$. If $x_{mn}^t=1$, then constraint \eqref{seq:hla} guarantees that $d_{mn}^t\geq 0$; constraint \eqref{seq:hlb} becomes an equality; and \eqref{seq:hlc} holds trivially. If $x_{mn}^t=0$, the flow changes direction $d_{mn}^t\leq 0$; constraint \eqref{seq:hlc} becomes an equality; and \eqref{seq:hlb} holds trivially. Observe that for $d_{mn}^t=0$, the indicator variable $x_{mn}^t$ becomes \emph{inconsequential}, and $h_m^t=h_n^t$ for any value of $x_{mn}^t$.

To maintain nodal pressures at desirable levels, water utilities use pumps installed on designated pipes to raise pressure. A water pipe equipped with a pump may be modeled as an ideal (lossless) pump followed by a pipe with pressure drop dictated by \eqref{eq:headlossMI}. The subset of edges representing ideal pumps is denoted by $\mcP_a\subset \mcP$. The remaining edges comprise the set $\bmcP_a:=\mcP\setminus\mcP_a$ and represent lossy pipes, for which the constraints in \eqref{eq:headlossMI} apply. Any reference to pump $(m,n)$ will henceforth refer to the ideal segment of the pump.

If pump $(m,n)\in \mcP_a$ is running during period $t$, its flow is constrained to lie within the range $\underline{d}_{mn}\leq d_{mn}^t \leq \overline{d}_{mn}$ with $\underline{d}_{mn}\geq 0$ due to engineering limitations \cite{Verleye_tankpumpmodel}. The pump $(m,n)$ adds pressure $g_{mn}^t\geq 0$ so that
\begin{equation}\label{eq:headgain}
h_n^t-h_m^t=g_{mn}^t.
\end{equation}
The pressure gain $g_{mn}^t$ depends on the pump speed and the water flow. This dependence is oftentimes approximated by a quadratic function \cite{Ulanicki_pumpmodel}, \cite{Verleye_tankpumpmodel}, \cite{TaylorOWFcones}. The dependence of $g_{mn}^t$ on water flow is relatively weak and may be ignored without significant loss of accuracy~ \cite{TaylorOWFcones}, \cite{CohenQHmodel}. Thus, for a fixed-speed pump, the pressure gain $g_{mn}$ is constant when the pump is running; and zero, otherwise. Oftentimes, when a pump is not running, water can flow freely in either directions through a bypass valve connected in parallel to the pump and without incurring any pressure difference \cite{CohenQHmodel}. The operation of a pump along with its bypass valve can be captured using the big-$M$ trick via the mixed-integer model for all $(m,n)\in \mcP_a$
\begin{subequations}\label{eq:pumpMI}
\begin{align}
h_m^t-h_n^t&=-g_{mn}x_{mn}^t\label{seq:pumpa}\\
-M(1-x_{mn}^t)&\leq d_{mn}^t-\tilde{d}_{mn}^t\leq M(1-x_{mn}^t)\label{seq:pumpc}\\
\underline{d}_{mn}x_{mn}^t&\leq \tilde{d}_{mn}^t \leq \overline{d}_{mn}x_{mn}^t\label{seq:pumpd}\\
x_{mn}^t&\in\{0,1\}.
\end{align}
\end{subequations}
The binary variable $x_{mn}^t$ indicates whether pump $(m,n)\in\mcP_a$ is running at time $t$. When the pump is running ($x_{mn}^t=1$), constraint \eqref{seq:pumpa} implies \eqref{eq:headgain}; otherwise ($x_{mn}^t=0$), it enforces $h_m^t=h_n^t$. For $x_{mn}^t=1$, constraints \eqref{seq:pumpc}--\eqref{seq:pumpd} imply that $\tilde{d}_{mn}^t=d_{mn}^t$ and the water flow in the pump is kept within the positive limits $[\underline{d}_{mn},\overline{d}_{mn}]$. For $x_{mn}^t=0$, variable $\tilde{d}_{mn}^t$ is set to zero and $d_{mn}^t$ represents the water flowing through the bypass valve of the pump. The auxiliary variable $\tilde{d}_{mn}^t$ will be useful later in computing the energy consumption of pump $(m,n)$. 

Note that a variable-speed pump model is not a generalization of a fixed-speed one unless non-trivial upper and lower bounds on the pump speeds are enforced. For instance, the OWF formulation for variable speed pumps in \cite{ZamzamOWPF,TaylorOWFcones} can not be used for fixed-speed pumps. Although there is an ongoing transition towards variable-speed pumps, the conventional WDS have a fleet of fixed-speed pumps which give way to on/off and implicit flow control \cite{odan2015epanet}, \cite{Verleye_tankpumpmodel}, \cite{Giacomello2013hybrid}. Thus, this work considers fixed-speed pumps.

The pressure at a reservoir can be assumed constant across days or weeks~\cite{TaylorOWFcones}. Consider reservoir $m\in\mcM_r$ whose constant pressure is $\bar{h}_m$. To draw water from this reservoir, its nodal pressure $h_m^t$ must be smaller than the constant pressure head $\bar{h}_m$ of the reservoir. This is enforced through the constraints
\begin{subequations}\label{eq:res}
	\begin{align}
	0\leq d_m^t&\leq M\alpha_m^t\label{seq:resA}\\
	h_m^t&\leq \bar{h}_m+ M(1-\alpha_m^t)\label{seq:resB}\\
	\alpha_m^t&\in\{0,1\}
	\end{align}
\end{subequations}
for all $m\in\mcM_r$ and times. The binary variable $\alpha_m^t$ indicates if water is drawn from reservoir $m$ at time $t$. If $\alpha_m^t=1$, reservoir $m$ is connected to the WDS and the constraints in \eqref{eq:res} ensure that $d_m^t\geq 0$ and $h_m^t\leq \bar{h}_m$. On the other hand, when $\alpha_m^t=0$, reservoir $m$ is disconnected, $d_m^t=0$, and constraint \eqref{seq:resB} is trivially satisfied.

As opposed to reservoirs, the water volume in tanks varies significantly during the day \cite{TaylorOWFcones}. Variations in water volume translate to variations in water level, which cause in turn variations in pressure at the bottom of the tank. To model the operation of tanks, let $\ell_m^t$ denote the water level in tank $m\in\mcM_b$ at the end of period $t$. To be consistent with the piezometric pressure head, the water level $\ell_m^t$ includes the geographical elevation of tank $m$. If $\delta$ is the duration of a control period and $A_m$ is the uniform cross-sectional area for tank $m$, the water level in tank $m$ satisfies the dynamics
\begin{equation}\label{eq:tank1}
\ell_m^t=\ell_m^{t-1}-\frac{d_m^t\delta}{A_m}.
\end{equation}
Due to its finite volume, the water level in tank $m$ is constrained at all times $t$ as
\begin{equation}\label{eq:tank2}
\underline{\ell}_m\leq \ell_m^t \leq \overline{\ell}_m.
\end{equation}
Typically, the net water exchange from tanks is kept at zero during the entire period of operation, that is 
\begin{equation}\label{eq:tank3}
\ell_m^0=\ell_m^T.
\end{equation} 

\begin{figure}[t]
	\centering
	\includegraphics[scale=0.36]{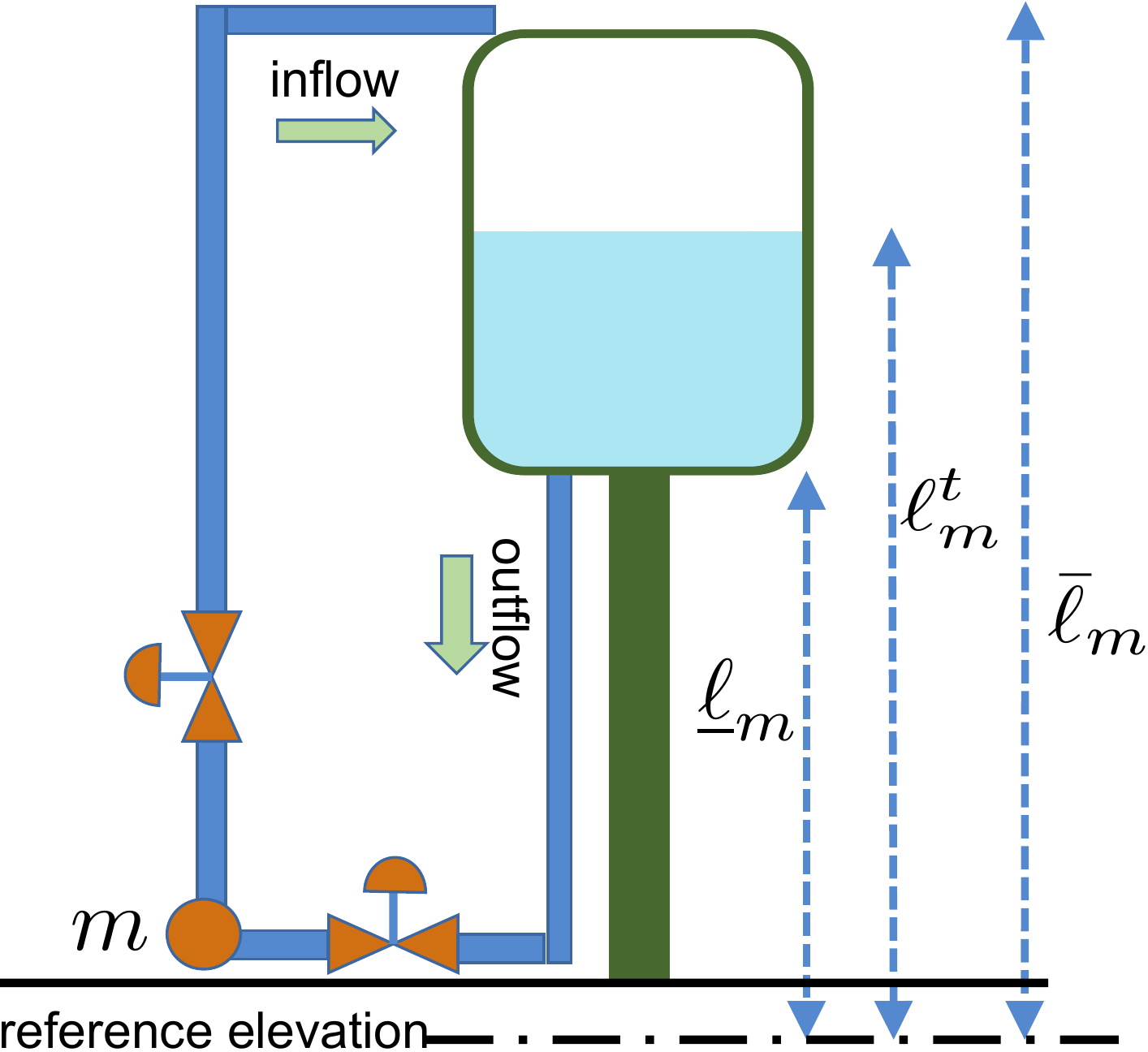}
	\caption{A schematic for a water tank sited at node $m$. The geographical elevation has been incorporated by referring heights to a common reference.}
	\label{fig:tank}
\end{figure}

Each tank has two separate paths for filling and emptying; see Fig.~\ref{fig:tank}. The filling or inlet pipe is connected near the top, and the emptying or outlet pipe is connected at the bottom. The two pipes are controlled by two separate valves. The output pressure of the valves can equal or less than the input pressure. Therefore, when tank $m$ is being filled in with water at time $t$, it should hold $h_m^t \geq \overline{\ell}_m$. Conversely, when water flows out of the tank, it follows that $h_m^t\leq \ell_m^t$. By closing both the inlet and outlet valves, the pressure $h_m^t$ at node $m$ becomes decoupled from the pressure at the bottom of the tank, $\ell_m^t$. 

To capture the aforementioned tank operation, let us introduce two binary variables $(\alpha_m^t,\beta_m^t)$ and the auxiliary continuous variable $\tilde{h}_m^t$. The operation of tank $m$ at time $t$ is described by the constraints
\begin{subequations}\label{tankMI}
	\begin{align}
	-M(1-\alpha_m^t)&\leq \tilde{h}_m^t-h_m^t \leq M(1-\alpha_m^t)\label{seq:tankMIA}\\
	-M\alpha_m^t&\leq d_m^t\leq M\alpha_m^t\label{seq:tankMIB}\\
	-M\beta_m^t&\leq d_m^t\leq M(1-\beta_m^t)\label{seq:tankMIC}\\
	\overline{\ell}_m- M(1-\beta_m^t)&\leq \tilde{h}_m^t\leq \ell_m^t+ M\beta_m^t\label{seq:tankMID}\\
\alpha_m^t,\beta_m^t&\in\{0,1\}.
	\end{align}
\end{subequations}
The variable $\alpha_m^t$ indicates if tank $m$ is connected at time $t$; and if it is, the variable $\beta_m^t$ indicates if the tank is filling. When the tank is connected  ($\alpha_m^t=1$), constraint \eqref{seq:tankMIA} yields $\tilde{h}_m^t=h_m^t$ and \eqref{seq:tankMIB} holds trivially. If additionally the tank is filling ($\beta_m^t=1$), then $d_m^t\leq 0$ from \eqref{seq:tankMIC} and $\tilde{h}_m^t=h_m^t\geq \overline{\ell}_m$ from \eqref{seq:tankMID}. If the tank is connected but emptying ($\alpha_m^t=1,\beta_m^t=0$), then $d_m^t\geq 0$ from \eqref{seq:tankMIC} and $\tilde{h}_m^t=h_m^t\leq \ell_m^t$ from \eqref{seq:tankMID}. When the tank is disconnected $(\alpha_m^t=0)$, constraint \eqref{seq:tankMIB} enforces $d_m^t=0$, the pressure in the tank is not related to the network pressure and the values of $\beta_m^t$ and $\tilde{h}_m^t$ are \emph{inconsequential}.

Valves are a vital flow-control component. Popular models for valves include an on/off switch model; a linear pressure-reducing model; and a flow-dependent nonlinear model~\cite{TaylorOWFcones}. Presuming a combination of on/off and linear valves on lossy pipes, a convex relaxation for OWF was put forth in \cite{TaylorOWFcones}. Although this simplistic setup can be incorporated here, this work addresses the more realistic WDS setup where valves are present only at reservoirs and tanks.

\section{Problem Formulation}\label{sec:PF}
With dynamic pricing, the objective here is to minimize the cost of electricity consumed by water pumps. This section collects the network constraints listed earlier and defines the OWF problem. The mechanical power consumed by pump $(m,n)\in\mcP_a$ during period $t$ in watts is given by the product of the induced pressure difference $g_{mn}$ measured in pascal, times the water flow $\tilde{d}_{mn}^t$ in m$^3$/sec~\cite{TaylorOWFcones}. If the overall energy efficiency of the pump is $\eta_{mn}$, it consumes electric energy $\frac{\delta\rho\tilde{g}g_{mn}}{\eta_{mn}}\tilde{d}_{mn}^t$ during time $t$ of duration $\delta$. For the fixed-speed pumps considered here, the pressure gain $g_{mn}$ is constant and we can thus define the electricity consumption coefficient \[c_{mn}:=\frac{\delta\rho\tilde{g}g_{mn}}{\eta_{mn}},\quad \forall (m,n)\in \mcP_a.\]

The OWF problem can be formally stated as follows. Given the initial water level in tanks $\{\ell_m^0\}_{m\in\mcM_b}$, the water demands at consumption nodes $\{d_m^t\}_{m\in\mcM\setminus\mcM_b\cup\mcM_r}$, the electricity prices $\{\pi_t\}_{t=1}^T$, and network parameters (tank capacities, pipe dimensions, pump pressure gains and minimum pressure requirements, tank heights); the OWF task aims at minimizing the electricity cost for running the pumps while meeting water demands and respecting WDS limitations. 

In detail, the pumping cost can be formulated as 
\begin{equation}
f(\tbd):=\sum_{t=1}^{T}\sum_{(m,n)\in \mcP_a} c_{mn}\pi_t\tilde{d}_{mn}^t
\end{equation}
where vector $\tbd$ collects the water flows $\{\tilde{d}_{mn}^t\}_{t}$ in all pumps $(m,n)\in\mcP_a$ and at all times. To simplify the presentation, the price of electricity $\pi_t$ is assumed invariant across the WDS for all $t$. The OWF problem can be posed as the minimization
\begin{align*}\label{eq:P1}
\min~&~f(\tbd)\tag{P1}\\
\mathrm{over}~& ~ {\{h_m^t}\}_{m\in\mcM}, \{d_m^t\}_{m \in \mcM_b\cup\mcM_r},\{d_{mn}^t\}_{(m,n)\in\mcP},\\
~&~\{\tilde{h}_m^t\}_{m\in\mcM_b},\{\ell_m^t\}_{m\in\mcM_b}, \{\tilde{d}_{mn}^t\}_{(m,n)\in\mcP_a},\\
~&~\{x_{mn}^t\}_{(m,n)\in\mcP},\{\alpha_m^t\}_{m\in\mcM_r\cup\mcM_b},\{\beta_m^t\}_{m\in\mcM_b},\quad  \forall t\\
\mathrm{s.to}~&~\eqref{eq:node1},\eqref{eq:hmin},\eqref{eq:headlossMI},\eqref{eq:pumpMI}-\eqref{tankMI}.
\end{align*}

Problem \eqref{eq:P1} involves the continuous variables $\{h_m^t,d_m^t,d_{mn}^t,\tilde{h}_m^t,\tilde{d}_{mn}^t\}$ and the binary variables $\{x_{mn}^t,\alpha_m^t,\beta_m^t\}$. For fixed-speed pumps, the cost in \eqref{eq:P1} is linear. Although most of the constraints are linear thanks to the big-$M$ trick, the constraints \eqref{seq:hlb}--\eqref{seq:hlc} modeling the pressure drop are non-linear. In fact, each one of these constraints involves one convex and one non-convex quadratic inequality. To obtain affordable OWF solutions, Section~\ref{sec:CR} relaxes the non-convex constraints and derives a mixed-integer problem that is convex with respect to the continuous variables.


\section{Convex Relaxation}\label{sec:CR}
The pressure drop across a lossy pipe $(m,n)\in\bmcP_a$ depends on its water flow $d_{mn}^t$ through the quadratic law of \eqref{eq:headloss}, which can be relaxed to a convex inequality as
\begin{itemize}
\item $h_m^t-h_n^t\geq c_{mn}(d_{mn}^t)^2$ for $d_{mn}^t\geq 0$; or
\item $h_n^t-h_m^t\geq c_{mn}(d_{mn}^t)^2$ for $d_{mn}^t\leq 0$. 
\end{itemize}
Since the sign of $d_{mn}^t$ is captured by the binary variable $x_{mn}^t$, the relaxation can be alternatively performed on \eqref{eq:headlossMI} to yield
\begin{subequations}\label{eq:relaxed}
	\begin{align}
	-&M(1-x_{mn}^t)\leq d_{mn}^t\leq Mx_{mn}^t\label{seq:hlra}\\
	-&M(1-x_{mn}^t)\leq h_m^t-h_n^t-c_{mn}(d_{mn}^t)^2\label{seq:hlrb}\\
	&h_m^t-h_n^t + c_{mn}(d_{mn}^t)^2\leq Mx_{mn}^t.\label{seq:hlrc}
	\end{align}
\end{subequations}
Comparing \eqref{eq:headlossMI} to \eqref{eq:relaxed}, the rightmost inequality of \eqref{seq:hlb} and the leftmost inequality of \eqref{seq:hlc} have been dropped in \eqref{eq:relaxed}. These are exactly the non-convex constraints. Replacing \eqref{eq:headlossMI} by \eqref{eq:relaxed} in \eqref{eq:P1}, leads to the relaxed problem 
\begin{align*}\label{eq:P2}
\min~&~f(\tbd)\tag{P2}\\
\mathrm{over}~& ~ {\{h_m^t}\}_{m\in\mcM}, \{d_m^t\}_{m \in \mcM_b\cup\mcM_r},\{d_{mn}^t\}_{(m,n)\in\mcP},\\
~&~\{\tilde{h}_m^t\}_{m\in\mcM_b},\{\tilde{d}_{mn}^t\}_{(m,n)\in\mcP_a},\\
~&~\{x_{mn}^t\}_{(m,n)\in\mcP},\{\alpha_m^t\}_{m\in\mcM_r\cup\mcM_b},\{\beta_m^t\}_{m\in\mcM_b},\quad  \forall t\\
\mathrm{s.to}~&~ \eqref{eq:node1},\eqref{eq:hmin},\eqref{eq:pumpMI}-\eqref{tankMI},\eqref{eq:relaxed}.
\end{align*}
Problem \eqref{eq:P2} is convex with respect to the continuous variables, and it could be handled by existing mixed-integer off-the-shelf solvers. Being a relaxation, the optimal value of \eqref{eq:P2} serves as a lower bound for the optimal value of \eqref{eq:P1}. If a minimizer of \eqref{eq:P2} satisfies \eqref{seq:hlrb} or \eqref{seq:hlrc} with equality for all $(m,n)\in\bmcP_a$, the relaxation is deemed \emph{exact}. In this case, the minimizer of \eqref{eq:P2} coincides with the minimizer of \eqref{eq:P1}. Nonetheless, the relaxation is not necessarily exact. 

To study the feasible sets of \eqref{eq:P1} and \eqref{eq:P2}, let $\bh$ collect the nodal pressures $\{h_m^t\}_{m,t}$; vector $\bd$ the water flows $\{d_{mn}^t\}_t$ for all $(m,n)\in\mcP$; and $\tbd$ has been defined after \eqref{eq:P1}. Define the projection of the feasible set of \eqref{eq:P1} into $(\tbd,\bd,\bh)$ as $\mcS_1$, and the projection of the feasible set of \eqref{eq:P2} into $(\tbd,\bd,\bh)$ as $\mcS_2$. The next result shows there exists a bijection between $\mcS_1$ [resp. $\mcS_2$] and the feasible set of \eqref{eq:P1} [resp. \eqref{eq:P2}].

\begin{lemma}\label{le:mapping1}
The $\bs:=\{\tbd,\bd,\bh\}$ components of any feasible point of \eqref{eq:P1} and \eqref{eq:P2} are sufficient to characterize the feasible point, modulo some inconsequential variables. 
\end{lemma}

\begin{IEEEproof}
It will be shown that upon fixing $(\tbd,\bd,\bh)$, the remaining variables listed under \eqref{eq:P1}--\eqref{eq:P2} can be determined, with only possible ambiguities on the values of inconsequential variables as detailed below. Given $\bd$, the water injections $\{d_n^t\}_{n,t}$ are set by \eqref{eq:node1}. Subsequently, the water levels $\{\ell_m^t\}_{m,t}$ are set by iterative computation of \eqref{eq:tank1} starting from the known initial tank level $\ell_m^0$. 

The binary variables capturing flow directions in lossy pipes can be recovered as
\[x_{mn}^t = \left\lfloor \frac{\sign(d_{mn}^t)+1}{2}\right\rfloor,\quad \forall (m,n)\in\bmcP_a,t\]
where $\lfloor a \rfloor$ denotes the floor function. If $d_{mn}^t=0$, the value of $x_{mn}^t$ is inconsequential and the aforementioned mapping sets it to zero. The binary variables pump statuses are set as $x_{mn}^t=\sign(\tilde{d}_{mn}^t)$ for $(m,n)\in\mcP_a$.

The variables governing reservoirs and tanks are set as 
\begin{subequations}\label{eq:le:r+t}
\begin{align}
\alpha_m^t&=|\sign(d_m^t)|,\quad \forall m\in\mcM_b\\
\beta_m^t&=\left\lfloor\frac{1-\sign(d_m^t)}{2}\right\rfloor,\quad \forall m\in\mcM_b\\
\tilde{h}_m^t&=\alpha_m^t h_m^t,\quad \forall m\in\mcM_b.
\end{align}
\end{subequations}
If tank $m$ is disconnected at time $t$, then $\alpha_m^t=0$ and the values of $\beta_m^t$ and $\tilde{h}_m^t$ become inconsequential. In that case, the mapping in \eqref{eq:le:r+t}  sets them to zero without harming feasibility.
\end{IEEEproof}

Lemma~\ref{le:mapping1} asserts that \eqref{eq:P1} and \eqref{eq:P2} can be equivalently expressed only in terms of $\bs:=\{\tbd,\bd,\bh\}$. The remaining variables have been introduced merely to avoid discontinuous or non-differentiable functions (e.g., sign or absolute value) as well as products between continuous and binary variables. In light of Lemma~\ref{le:mapping1} and with a slight abuse in terminology, we will henceforth refer to $\mcS_1$ [resp. $\mcS_2$] as the \emph{feasible set} of \eqref{eq:P1} [resp. \eqref{eq:P2}]. Due to the relaxation, it holds $\mcS_1\subseteq \mcS_2$.

When it comes to \eqref{eq:P1}, a feasible point can be constructed only by its $\{\tbd,\bd\}$ components, since a feasible $\bh$ can be recovered from $\{\tbd,\bd\}$ as follows. Given $\{\tbd,\bd\}$, the variables $\{x_{mn}^t,\alpha_{m}^t,\beta_{m}^t,d_m^t,\ell_m^t\}$ can be set as in the proof of Lemma~\ref{le:mapping1}. The values of \emph{pressure differences} across pipes can be found by \eqref{eq:headlossMI} and \eqref{seq:pumpa}. The next question is how to recover pressures from pressure differences. 

To express pressure differences at time $t=1,\ldots,T$, let us define an edge-node incidence matrix depending on the water flow directions at time $t$. Define $\bd^t$ as the subvector of $\bd$ collecting water flows only at time $t$. Then, introduce the $P\times |\mcM|$ incidence matrix $\bA(\bd^t)$ so that if its $p$-th row corresponds to pipe $p=(m,n)$, then its $(p,k)$ entry is
\begin{equation*}
A_{p,k}(\bd^t):=\left\{\begin{array}{ll}
-\sign^2(d_{mn}^t)+\sign(d_{mn}^t)+1	&,~k=m\\
\sign^2(d_{mn}^t)-\sign(d_{mn}^t)-1		&,~k=n\\
0												&,~\textrm{otherwise}.
\end{array}\right.
\end{equation*}
In this way, vector $\bA(\bd^t)\bh^t$ captures the pressure differences taken across the direction of water flows. For zero flows, the standard pipe direction $(m,n)$ is selected without loss of generality.

If $(\bh^t,\tbd^t)$ are the subvectors of $(\bh,\tbd)$ corresponding to time $t$, the pressure differences can be expressed as
\begin{equation}\label{eq:heq}
\bA(\bd^t)\bh^t=\bb(\tbd^t,\bd^t),\quad \forall t
\end{equation}
where $\bb(\tbd^t,\bd^t)$ is the mapping induced by \eqref{eq:headlossMI} and \eqref{seq:pumpa}. Since $\{\tbd,\bd\}$ is feasible for \eqref{eq:P1}, the overdetermined system in \eqref{eq:heq} is consistent. However, its solution is not unique: The all-one vector $\mathbf{1}$ belongs to the nullspace of $\bA(\bd^t)$ by definition, so if $\bh^t$ satisfies \eqref{eq:heq}, then $\bh^t+c\mathbf{1}$ satisfies \eqref{eq:heq} too for any $c$.

Satisfying \eqref{eq:heq} alone is not sufficient for $\bh^t$ to be feasible for \eqref{eq:P1}. It should also satisfy the inequality constraints \eqref{eq:hmin}, \eqref{seq:resB}, \eqref{seq:tankMIA}, and \eqref{seq:tankMID}. These constraints are abstractly expressed as
\begin{equation}\label{eq:hineq}
\underline{\bh}(\tbd,\bd)\leq \bh\leq \overline{\bh}(\tbd,\bd).
\end{equation}


Given $\{\tbd,\bd\}$ for a feasible point of \eqref{eq:P1}, a feasible pressure vector $\bh$ can be found by ensuring \eqref{eq:heq}--\eqref{eq:hineq}. A water utility would implement $\bh$ by controlling the pressures at reservoir valves. The aforesaid procedure proves the following claim.

\begin{lemma}\label{le:mapping2}
	Any feasible point of \eqref{eq:P1} is characterized by its $\{\tbd,\bd\}$ components modulo some inconsequential variables. A vector of feasible pressures $\bh$ can be recovered by solving the linear program (LP)
	\begin{align}\label{eq:LP}
    \mathrm{find}~&~\bh\\
     \mathrm{s.to}~&~\eqref{eq:heq}-\eqref{eq:hineq}.\nonumber
	\end{align}
\end{lemma}

Let $\mcH(\tbd,\bd)$ be the set of vectors $\bh$ solving the feasibility problem in \eqref{eq:LP}. Lemma \ref{le:mapping2} implies that any solution to \eqref{eq:LP} provides a feasible point for \eqref{eq:P1}.

\begin{figure}[t]
	\centering
	\includegraphics[scale=0.4]{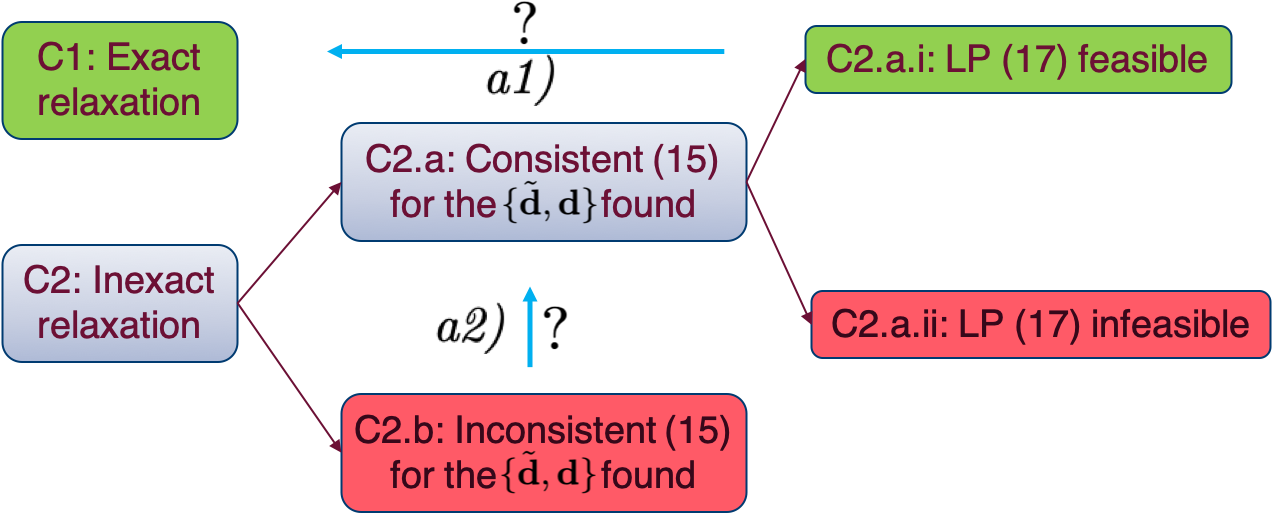}
	\caption{Possible cases for the feasibility of a minimizer obtained by \eqref{eq:P2}. Problem \eqref{eq:P3} converts case \emph{C2.a.i} to \emph{C1}. Moreover, under the conditions of Lemma~\ref{le:parallel}, it also converts case \emph{C2.b} to \emph{C2.a}.}
	\label{fig:cases}
\end{figure}

Given Lemma~\ref{le:mapping2}, let us see if one can find a feasible point for \eqref{eq:P1} by solving \eqref{eq:P2}. Consider a minimizer $\bs_1:=\{\tbd_1,\bd_1,\bh_1\}$ of \eqref{eq:P1} attaining the cost $f_1:=f(\tbd_1)$. Consider also a minimizer $\bs_2:=\{\tbd_2,\bd_2,\bh_2\}$ of \eqref{eq:P2} with $f_2:=f(\tbd_2)$ with $f_2\leq f_1$ due to the relaxation. The next cases can be identified for $\bs_2$ as illustrated in Figure~\ref{fig:cases}:
\renewcommand{\labelenumi}{\textit{C\arabic{enumi}}.}
\renewcommand{\labelenumii}{\textit{\labelenumi\alph{enumii}}.}
\renewcommand{\labelenumiii}{\textit{\labelenumii\roman{enumiii}}.}
\begin{enumerate}
	\item If the relaxation is exact, then $\bh_2\in\mcH(\tbd_2,\bd_2)$; the costs agree $f_2= f_1$; and $\bs_2$ can be implemented in lieu of $\bs_1$.
	\item If the relaxation is inexact, vector $\bh_2$ satisfies only the equations in \eqref{eq:heq} related to pumps, whereas some of the constraints related to lossy pipes in \eqref{eq:relaxed} are satisfied with strict inequalities. In this case, one may try to recover a vector of physically feasible pressures by enforcing \eqref{eq:heq}--\eqref{eq:hineq}. The following subcases are identified.
		\begin{enumerate}
		\item The linear system of \eqref{eq:heq} is \emph{consistent} for $\bb(\tbd_2,\bd_2)$. Again, two cases can be identified.
		\begin{enumerate}
			\item The LP in \eqref{eq:LP} is feasible for $(\tbd_2,\bd_2)$ with $\cbh_2\in\mcH(\tbd_2,\bd_2)$. The point $\cbs_2:=\{\tbd_2,\bd_2,\cbh_2\}$ is feasible for \eqref{eq:P1} and attains the cost $\check{f}_2:=f(\tbd_2)=f_2$. Because $\cbs_2$ is feasible for \eqref{eq:P1}, the optimal cost has been attained, that is $\check{f}_2=f_2=f_1$.
			\item The LP in \eqref{eq:LP} is infeasible for $(\tbd_2,\bd_2)$. A feasible point for \eqref{eq:P1} cannot be recovered.
		\end{enumerate}
		\item The linear system of \eqref{eq:heq} is \emph{inconsistent} for $\bb(\tbd_2,\bd_2)$. A feasible point for \eqref{eq:P1} cannot be recovered.
	\end{enumerate}	
\end{enumerate}

Cases \textit{C1} and \textit{C2.a.i} are computationally useful since they recover an optimal point. On the other hand, cases \textit{C2.a.ii} and \textit{C2.b}, do not provide any useful output. Based on numerical tests with different WDS networks and under various pricing/demand scenarios, we have empirically observed that:
\begin{itemize}
\item Case \emph{C1} occurs rarely.
\item Case \emph{C2.a.i} is encountered frequently in radial networks.
\item Case \emph{C2.b} occurs frequently in meshed networks. 
\end{itemize}
Spurred by these observations and to improve the chances for an exact relaxation of \eqref{eq:P1}, the next section adds a penalty term in the objective of \eqref{eq:P2}. It then studies the feasibility and optimality of this penalized convex relaxation.

\section{Penalized Convex Relaxation}\label{sec:PCR}
Toward an exact relaxation of \eqref{eq:P1}, define the penalty 
\begin{equation}\label{eq:penalty}
g(\bh):=\sum_{t=1}^{T}\sum_{(m,n)\in \bmcP_a}|h_m^t-h_n^t|
\end{equation}
which sums up the absolute pressure differences across lossy pipes and over all times. Let us formulate a \emph{penalized convex relaxation} by replacing the cost of \eqref{eq:P2} by
\begin{align*}\label{eq:P3}
\min~&~f(\tbd)+\lambda g(\bh)\tag{P3}\\
\mathrm{s.to}~&~\eqref{eq:node1},\eqref{eq:hmin},\eqref{eq:pumpMI}-\eqref{tankMI},\eqref{eq:relaxed}
\end{align*}
for $\lambda>0$. Sections~\ref{subsec:feas} and \ref{subsec:opt} next study respectively the feasibility and optimality of \eqref{eq:P3}.

\subsection{Improving Feasibility}\label{subsec:feas}
Although \eqref{eq:P2} and \eqref{eq:P3} share the same feasible set, this section shows that \eqref{eq:P3} features two advantages over \eqref{eq:P2} as depicted in Figure~\ref{fig:cases}:
\renewcommand{\labelenumi}{\textit{a\arabic{enumi}})}
\begin{enumerate}
\item Problem~\eqref{eq:P3} eliminates the occurrence of \emph{C2.a-i}. The problem instances falling under \emph{C2.a-i} with \eqref{eq:P2}, fall under the useful case \emph{C1} for \eqref{eq:P3}.
\item Under some conditions, problem \eqref{eq:P3} does not encounter the unfavorable case \emph{C2.b} either. 
\end{enumerate}

The following result establishes advantage \emph{a1)} and is shown in the appendix.

\begin{theorem}\label{th:1}
If $\bs_3:=\{\tbd_3,\bd_3,\bh_3\}$ is a minimizer of \eqref{eq:P3} and $\mcH(\tbd_3,\bd_3)$ is non-empty, then  $\bh_3\in\mcH(\tbd_3,\bd_3)$.
\end{theorem}

From Theorem~\ref{th:1} and Lemma~\ref{le:mapping2}, the next result follows.

\begin{corollary}\label{cor:th1}
Under the assumptions of Theorem~\ref{th:1}, the minimizer $\bs_3:=\{\tbd_3,\bd_3,\bh_3\}$ of \eqref{eq:P3} is feasible for \eqref{eq:P1}.
\end{corollary}

Corollary~\ref{cor:th1} asserts that if the water flows obtained from \eqref{eq:P3} can be mapped to physically feasible pressures, then the minimizer of \eqref{eq:P3} contains already physically feasible pressures and this shows advantage \emph{a1)}. In other words, instead of having to solve \eqref{eq:P2} first and then \eqref{eq:LP} to recover a feasible OWF schedule, a feasible schedule can be found by solving \eqref{eq:P3} alone.

Before moving to \emph{a2)}, some graph theory preliminaries are reviewed. Given an undirected graph $\mcG:=(\mcM, \mcP)$, its \emph{degree} is the number of incident edges. A graph is connected if there exists a sequence of adjacent edges between any two of its nodes. A minimal set of edges $\mcP_\mcT$ preserving the connectivity of a connected graph constitutes a \emph{spanning tree} of $\mcG$; is denoted by $\mcT:=(\mcM,\mcP_\mcT)$; and apparently, $|\mcP_\mcT|=|\mcM|-1$. The edges not belonging to a spanning tree $\mcT$ are referred to as \emph{links} with respect to $\mcT$. A \emph{cycle} is a sequence of adjacent edges without repetition that starts and begins at the same node. A \emph{tree} is a connected graph with no cycles. In a directed graph, each edge is assigned a directionality. A \emph{path} from node $m$ to $n$ is defined as a sequence of directed edges originating from $m$ and terminating at $n$. Given the undirected graph $(\mcM,\mcP)$ modeling a WDS and the vector $\bd^t$ of flows at time $t$, let us define the \emph{directed} graph $(\mcM,\mcP(\bd^t))$ where edge $p$ runs from node $m$ to node $n$ if $d_{m,n}^t\geq0$; and vice versa, otherwise.

To show \emph{a2)}, we study the consistency of \eqref{eq:heq}. Had the WDS graph been a tree, the edge-node incidence matrix would have been full row-rank~\cite{GodsilRoyle}. Hence, the equations in \eqref{eq:heq} would have been consistent for any $\bb(\tbd^t,\bd^t$). This implies that possible inconsistencies in \eqref{eq:heq} arise from cycles in $\mcG$. Because studying the generic case of cycles is not obvious, we consider the special case of a cycle where all but one nodes have degree two. This subset of edges will be henceforth termed a \emph{ring}. A ring can be rooted at the node with degree larger than two. We provide conditions under which a minimizer of \eqref{eq:P3} satisfies the constraints in \eqref{eq:relaxed} with equality for all edges of a ring.

\begin{lemma}\label{le:parallel}
Let $\bs_3=\{\tbd_3,\bd_3,\bh_3\}$ be a minimizer of \eqref{eq:P3} and $\bd_3^t$ be the subvector of $\bd_3$ collecting the flows at time $t$. If the directed graph $(\mcM,\mcP(\bd_3^t))$ contains a ring $\mcR\subseteq \mcP(\bd_3^t)$ rooted at node $m$, such that
\begin{itemize}
	\item all nodes incident to $\mcR$ have identical pressure limit $\underline{h}$;
	\item all nodes incident to $\mcR$ but $m$ host no tanks or reservoirs;
	\item all edges in $\mcR$ host no pumps;
\end{itemize} 
then $h_i^t-h_j^t=c_{ij}(d_{ij}^t)^2$ for all directed edges $(i,j)$ in $\mcR$.
\end{lemma}

Leveraging Lemma~\ref{le:parallel}, the ensuing result shows the advantage \emph{a2)} of \eqref{eq:P3} over \eqref{eq:P2} for a large class of WDS.

\begin{theorem}\label{th:2}
Let $\bs_3:=\{\tbd_3,\bd_3,\bh_3\}$ be a minimizer of \eqref{eq:P3} and $(\tbd_3^{t},\bd_3^{t})$ be the subvectors of $(\tbd_3,\bd_3)$ corresponding to time $t$. The system of equations in \eqref{eq:heq} is consistent for $\bs_3$ at time $t$, if all undirected cycles in $\left(\mcM,\mcP(\bd_3^t)\right)$ constitute rings satisfying the conditions of Lemma~\ref{le:parallel}.
\end{theorem}

To appreciate the claim of Theorem~\ref{th:2}, recall that for a point to be feasible for \eqref{eq:P1}, it is sufficient to satisfy \eqref{eq:heq} and \eqref{eq:hineq}. Since $\bA(\bd^t)\boldsymbol{1}=\mathbf{0}$, the next result can be inferred.

\begin{corollary}\label{cor:th2}
Under the assumptions of Theorem~\ref{th:2}, if the left or right inequality in \eqref{eq:hineq} are omitted, then a minimizer of \eqref{eq:P3} is feasible for \eqref{eq:P1}.
\end{corollary}

Corollary~\ref{cor:th2} asserts that \eqref{eq:P3} can be advantageous for coping with OWF tasks with no upper bounds on pressures; see also~\cite{TaylorOWFCDC}. An important problem complying to this setup is the \emph{water flow} (WF) task. Different from OWF, the WF problem solves the WDS equations over a single period upon specifying nodal water demands and a reference pressure. In a recent work~\cite{SinghKekatosWF19}, we have dealt with the WF task using a similar penalization, which is shown to yield the unique WF solution for a broader class of WDS.

\subsection{Optimality}\label{subsec:opt}
The previous section documented the advantages of \eqref{eq:P3} over \eqref{eq:P2} in terms of providing physically feasible OWF schedules under the conditions of Lemma~\ref{le:parallel} and Theorem~\ref{th:2}. However, the objective in \eqref{eq:P3} differs from the one of \eqref{eq:P1}: If a minimizer $\bs_3=\{\tbd_3,\bd_3,\bh_3\}$ of \eqref{eq:P3} is feasible for \eqref{eq:P1}, it will achieve in general a larger pumping cost than a minimizer of \eqref{eq:P1}, that is $f(\tbd_3)\geq f_1$. However, this suboptimality gap diminishes for decreasing $\lambda$ as explained next. We first review a general result on bi-objective optimization~\cite[Sec.~4.7.5]{BoVa04}:

\begin{lemma}[\cite{BoVa04}]\label{le:lambda}
Consider the minimization problem
	\[\bx_\lambda:=\arg\min_{\bx\in\mcX} f_a(\bx)+\lambda f_b(\bx),\]
for some real valued functions $f_a(\bx)$ and $f_b(\bx)$ defined on $\mcX$. If $\lambda_2>\lambda_1\geq 0$, then $f_a(\bx_{\lambda_2})\geq f_a(\bx_{\lambda_1})$.
\end{lemma}


Identifying functions $(f_a,f_b)$ of Lemma~\ref{le:lambda} to functions $(f,h)$ in the objective of \eqref{eq:P3} implies that for decreasing $\lambda$, a minimizer of \eqref{eq:P3} gives lower $f(\tbd_3(\lambda))$. However, the feasibility of $\bs_3$ for \eqref{eq:P1} is not guaranteed. If the conditions of Lemma~\ref{le:parallel} and Theorem~\ref{th:2} are met and $\bs_3$ is feasible for \eqref{eq:P1}, then $f(\tbd_3)\geq f_1$. Next, for $\lambda=0$, problem \eqref{eq:P3} degenerates to \eqref{eq:P2}, and gives a lower bound on $f_1$. Overall, we get that
\begin{equation}\label{eq:bottomline}
f(\tbd_2)\leq f_1\leq f(\tbd_3(\lambda)).
\end{equation}
From Theorems \ref{th:1} and \ref{th:2}, the advantage of the penalty term $g(\bh)$ does not depend on the value of $\lambda$ as long as $\lambda>0$. So under the conditions of Lemma~\ref{le:parallel} and Theorem~\ref{th:2}, one can choose arbitrarily small $\lambda$ to tighten the right-hand inequality in \eqref{eq:bottomline}. The caveat behind the bounds of \eqref{eq:bottomline} are the conditions assumed by Lemma~\ref{le:parallel} and Theorem~\ref{th:2}. Even though these conditions were grossly violated during the tests of Section~\ref{sec:tests}, the inequalities in \eqref{eq:bottomline} were frequently tightened to equalities. Albeit \eqref{eq:P2} oftentimes attained the optimal cost $f_1$, its minimizer was not feasible for \eqref{eq:P1}. In fact, there is no obvious way of converting the minimizer of \eqref{eq:P2} to a feasible point. Instead, problem \eqref{eq:P3} found a minimizer for \eqref{eq:P1} in most of the tests.

\section{Numerical Tests}\label{sec:tests} 

\begin{figure}[t]
	\centering
	\includegraphics[scale=0.28,angle=-90]{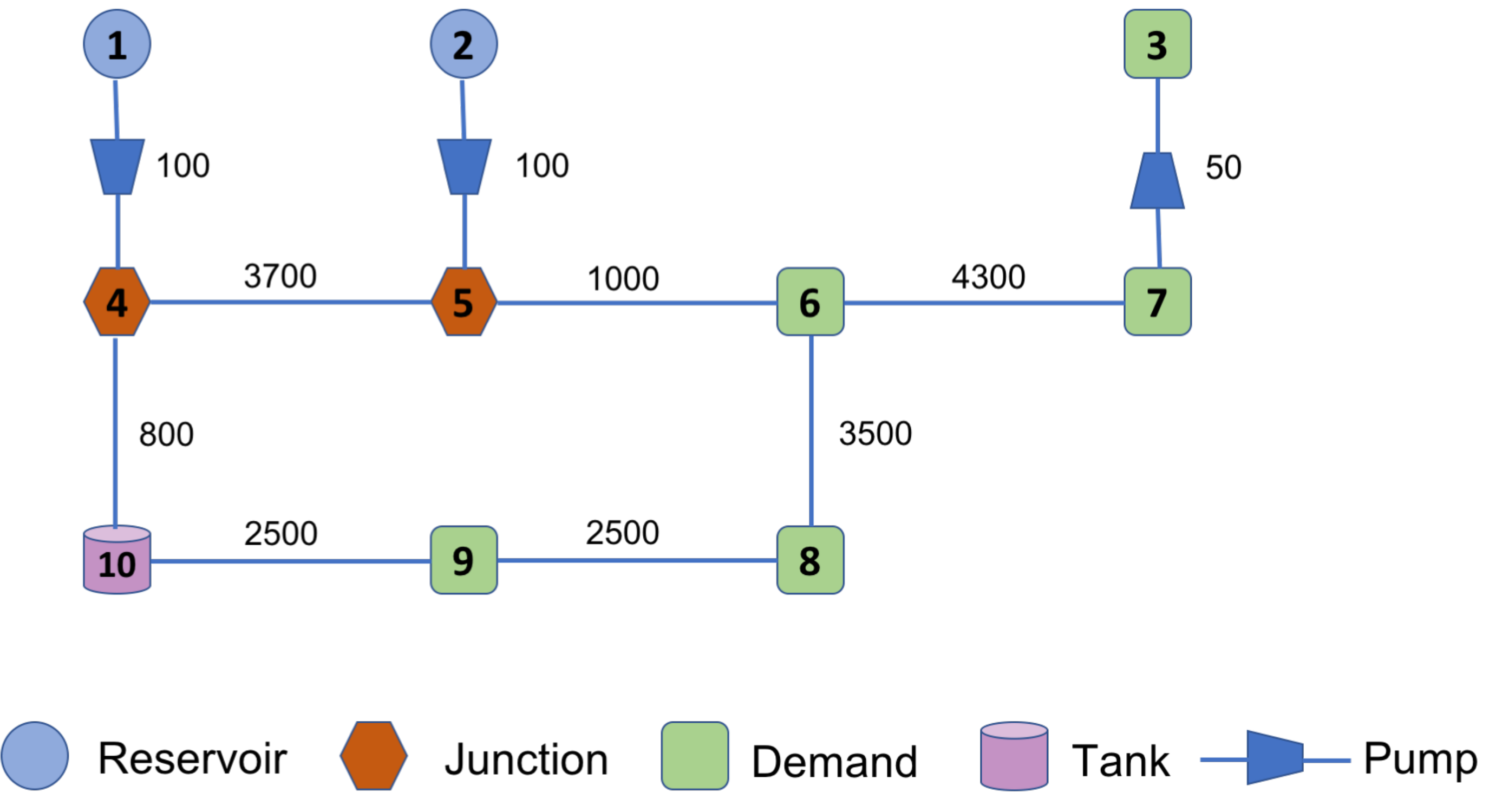}
	\caption{Benchmark water distribution system. The length for lossy pipes and head gain for pumps are shown in meters.}
	\label{fig:wds}
\end{figure}

\begin{figure}[t]
	\centering
	\includegraphics[scale=0.21]{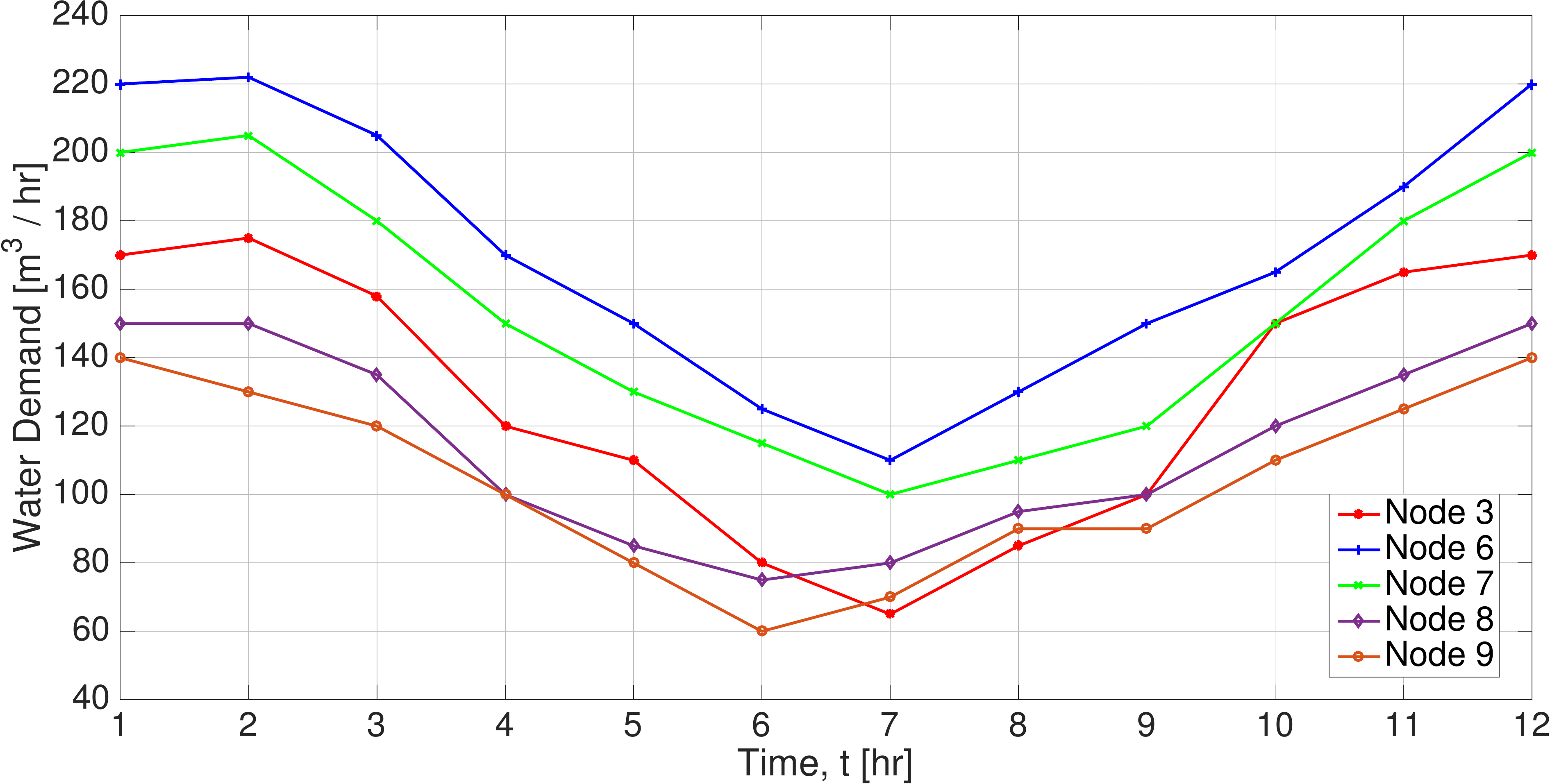}
	\caption{Per-node water demand across time.}
	\label{fig:demand}
\end{figure}

\begin{figure}[t]
	\centering
	\includegraphics[scale=0.22]{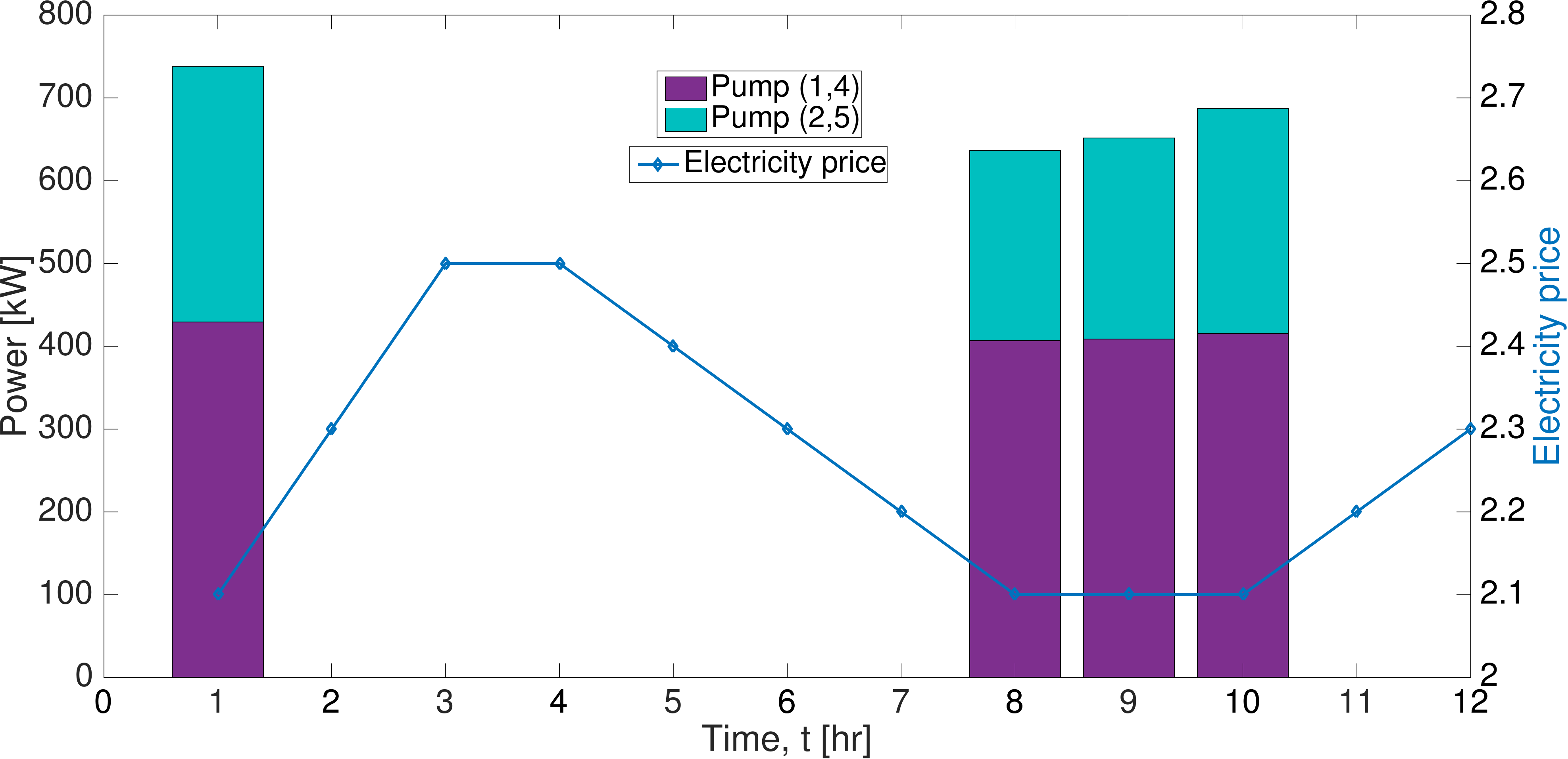}
	\includegraphics[scale=0.22]{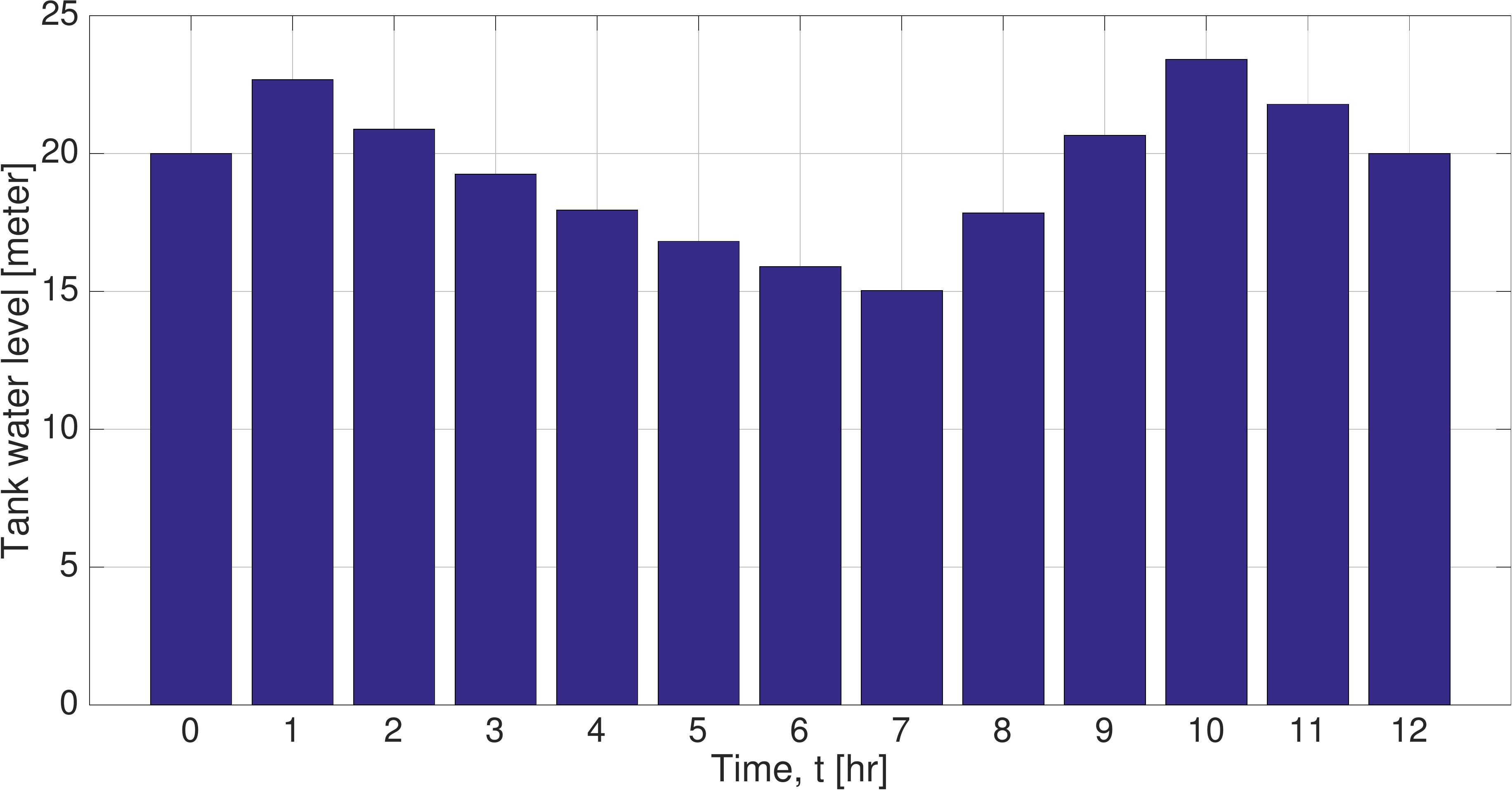}
	\caption{\emph{Top:} Electric power consumed by pumps during hour $t$. Pumps $(1,4)$ and $(2,5)$ were turned on during the same hours of lower electricity prices, whereas pump $(3,7)$ was not operated. Albeit the two pumps add the same pressure gain, they exhibit different electricity consumption due to different water flows. \emph{Bottom}: Water level in tank node $10$ at the end of hour $t$.}
	\label{fig:pump+tank}
\end{figure}

The new OWF solver was evaluated on the benchmark WDS of \cite{ZamzamOWPF}, \cite{CohenQHmodel}, which is shown in Figure~\ref{fig:wds}. It consists of $10$ nodes including $2$ reservoirs and a tank; $3$ fixed-speed pumps; and $7$ lossy pipes. All lossy pipes have a diameter of $0.4$m and friction coefficient $f_{m,n}=0.01$. The efficiency for all pumps is $85\%$ and for their motors $95\%$, resulting in an overall efficiency of $\eta=0.81$. The minimum and maximum water flows for all pumps are $100$m$^3/$hr and $1,500$m$^3/$hr, respectively. The pressure at reservoir nodes 1 and 2 is accordingly $-2.5$m and $5$m. The minimum pressure requirement $\underline{h}_m$ for nodes $3$ to $10$ is $\{10,7,12,10,5,10,10,10\}$m. Tank node $10$ has an area of $A_{10}=490.87\text{m}^2$; water level limits $\underline{\ell}_{10}=10$ and $\overline{\ell}_{10}=30$m; and initial water level $\ell_{10}^0=20$m.

The WDS was scheduled hourly for a horizon of $T=12$ hours for the demands of  Figure~\ref{fig:demand}; see~\cite{ZamzamOWPF}. The prices $\{\pi_t\}_{t=1}^{12}$ were set to the average day-ahead locational marginal prices during 8:00--20:00 on April 1,~2018 from the PJM market, and are shown in Fig.~\ref{fig:pump+tank}. The OWF tests were solved using the MATLAB-based optimization toolbox YALMIP along with the mixed-integer solver Gurobi~\cite{yalmip}, \cite{gurobi}. All tests were run on a $2.7$ GHz, Intel Core i5 computer with $8$ GB RAM.

We first checked whether the convex relaxation was exact. A minimizer of \eqref{eq:P3} was deemed feasible for \eqref{eq:P1} if ${|h_m^t-h_n^t|-c_{mn}\left(d_{mn}^t\right)^2\leq10^{-4}}$ for all pipes and times. A minimizer for \eqref{eq:P3} was obtained in $8.34$ sec for $\lambda=0.1$. The minimizer was in fact feasible for \eqref{eq:P1}. Figure~\ref{fig:pump+tank} presents the power consumed by pumps (top) and the water level in tank $10$ (bottom). The pumps run for the hours with the lowest prices over which tank node $10$ is filled, as expected. The tank is emptied during the hours of higher electricity prices, and its level is brought to its initial level at the end of the horizon. 

The modeling accuracy of the minimizer obtained by \eqref{eq:P3} was also tested against the standard simulation software EPANET~\cite{epanet2000}.The water injections obtained for the previous example by our MI-SOCP-based solver were fed into the water flow solver of EPANET to calculate the related pressures over the standard network model. The pressures found by the two models differed only by $0$--$0.91$ft across all nodes and times, with the median deviation being $0.21$ft. These differences are relatively insignificant considering that the average nodal pressure is on the order of $35$ft.

\begin{table}[t]
\renewcommand{\arraystretch}{1.4}
	\centering
	\caption{Pumping Cost Attained by \eqref{eq:P3} for Different $\lambda$'s}
	\vspace*{-1em}
	\begin{tabular}{|c| c| c|c|c|} 
		\hline\hline
		$\lambda$ & 0 & 0.01 & 0.1& 1  \\ \hline
		$f(\tbd_3)$ & 5,699.0& 5,699.0& 5,699.0& 5,704.2 \\ \hline
		\emph{comment} & lower bound \eqref{eq:P2} & infeasible &feasible&feasible \\ \hline
		\hline
	\end{tabular}
	\label{tbl:lambda}
\end{table} 

We next evaluated the effect of $\lambda$ on the feasibility and optimality of a minimizer of \eqref{eq:P3} with respect to \eqref{eq:P1}. We first solved \eqref{eq:P2} to obtain a lower bound $f(\tbd_2)$ on $f_1$. As a heuristic for setting $\lambda$, we computed $S:=\sum_{t=1}^{T}\sum_{(m,n)\in \bmcP_a}c_{mn}(d_{mn}^t)^2$ from the minimizer of \eqref{eq:P2}, and chose $\lambda=1$ so that $\lambda S$ was approximately $f(\tbd_2)/100$. For $\lambda=1$, the minimizer of \eqref{eq:P3} was feasible for \eqref{eq:P1} and provided an upper bound for $f_1$. To tighten \eqref{eq:bottomline}, problem \eqref{eq:P3} was solved for decreasing values of $\lambda$ obtaining the results of Table~\ref{tbl:lambda}. The minimizer of \eqref{eq:P3} for $\lambda=0.1$ was feasible for \eqref{eq:P1} and attained the same pumping cost as $f(\tbd_2)$. The infeasibility observed for $\lambda=0.01$ is attributed to the numerical accuracy of the solver, and such cases could be avoided by increasing $\lambda$. Hence, the minimizer of \eqref{eq:P3} constitutes a minimizer for \eqref{eq:P1} as well. It is worth stressing that even though the benchmark WDS of Figure~\ref{fig:wds} does not meet the conditions of Lemma~\ref{le:parallel} and Theorem~\ref{th:2}, an exact relaxation has been achieved. 

\begin{table*}[t]
	\renewcommand{\arraystretch}{1.0}
	\centering
	\caption{Suboptimality Gap Attained by Feasible Points Obtained through \eqref{eq:P3}}
	\vspace*{-1em}
	\begin{tabular}{|c| c| c|c|c|c|c|c|c|c|c|} 
		\hline
		Day of March 2018 & 10 & 11 & 12 & 13  & 14 & 15 & 16 & 17 & 18 & 19  \\ \hline
		$f(\tbd_2)$ & $6,968.5$ & $6,915.0$ & $8,524.6$ & $8,404.6$ & $8,220.5$ & $7,237.9$ & $7,206.8$ & $6,807.4$ & $6,404.0$ & $7,206.8$ \\ \hline
		$f(\tbd_3)$ & $7,042.8$ & $7,010.9$ & $8,524.6$ &$8,404.6$&$8,461.8$&$7,264.7$ & $7,206.8$ & $6,807.4$ & $6,527.1$ & $7,206.8$ \\ \hline
$\frac{f(\tbd_3)-f(\tbd_2)}{f(\tbd_2)}~[\%]$	&	$1.06$ & $1.39$ & $7\cdot10^{-9}$ & $3\cdot10^{-7}$ & $2.93$ & $0.37$ & $2\cdot10^{-7}$ & $7.\cdot10^{-4}$ & $1.92$ & $2\cdot10^{-7}$\\ \hline
$\lambda$ & $5$ & $5$ & $0.5$ & $1$ & $10$ & $2$ & $0.2$ & $0.83$ & $6$ & $0.6$\\ \hline
		Solution time [min:sec] & $00:07$ & $21:00$ & $00:09$ & $00:06$ & $22:19$ & $00:29$ & $00:10$ & $00:51$ & $20:50$ & $00:18$\\
		\hline
	\end{tabular}
	\label{tbl:march}
	\vspace*{-1em}
\end{table*} 

\begin{figure}[t]
	\centering
	\includegraphics[scale=0.23]{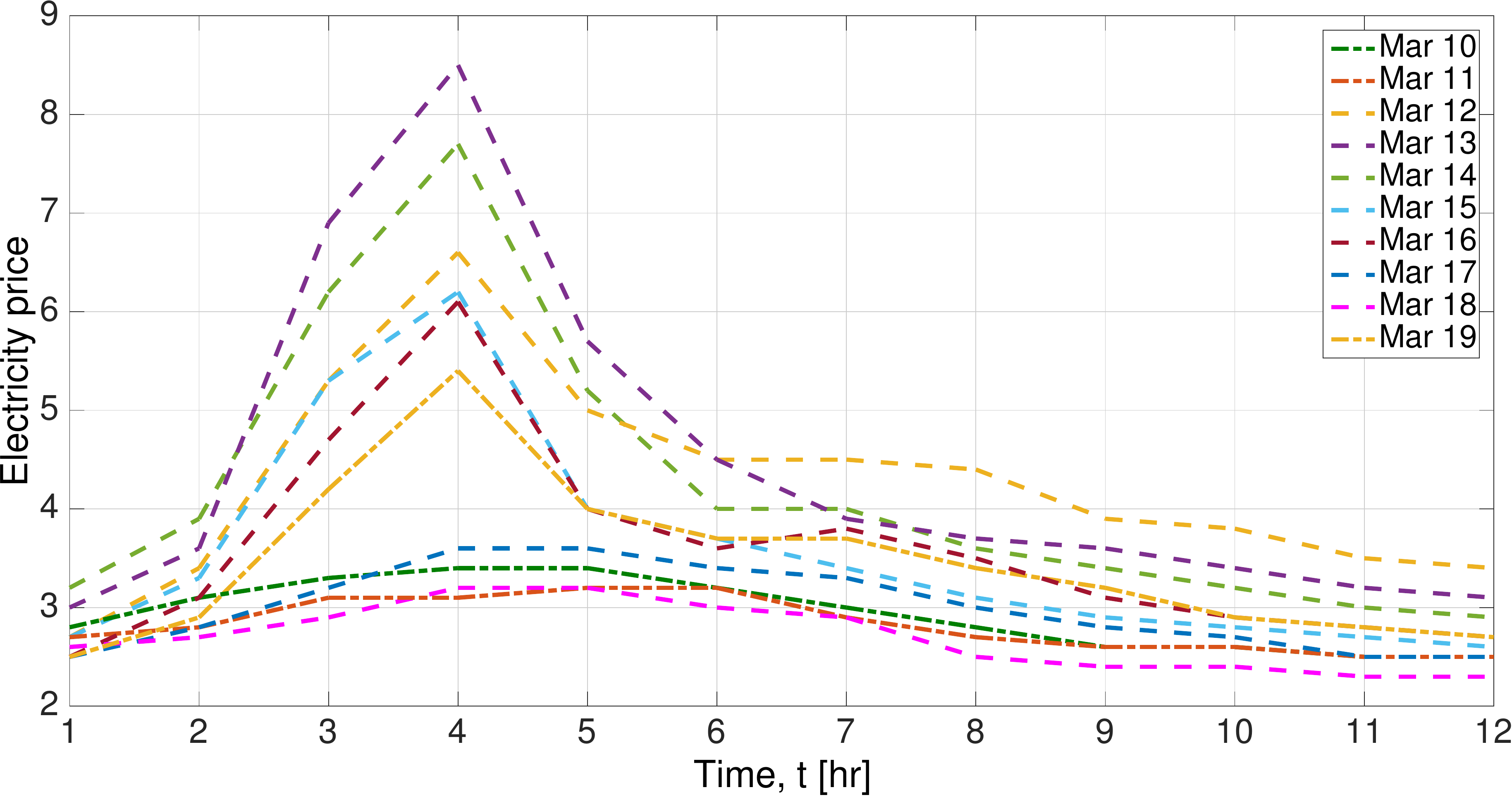}
	\caption{Day-ahead PJM electricity prices [$\cent$/kWh] for March 10--19, 2018.}
	\label{fig:prices}
\end{figure}

Similar tests were conducted for the PJM prices between March 10--19, 2018 during 5:00--17:00 shown in Fig.~\ref{fig:prices}. The results are summarized in Table~\ref{tbl:march}. For all $10$ days, problem \eqref{eq:P3} succeeded in finding a feasible point for the values of $\lambda$ reported in Table~\ref{tbl:march}. Moreover, the upper and lower bounds $f(\tbd_3)$ and $f(\tbd_2)$ were close implying small suboptimality gaps. It is worth stressing that the relaxation in \eqref{eq:P2} was \emph{inexact} for all tests. Albeit cost $f(\tbd_2)$ was equal to $f(\tbd_3)$ (and therefore equal to the optimal cost $f_1$ as well) for some cases, there is no obvious way to obtain an OWF dispatch from the minimizer of \eqref{eq:P2}.

\begin{figure}[t]
	\centering
	\includegraphics[scale=0.4]{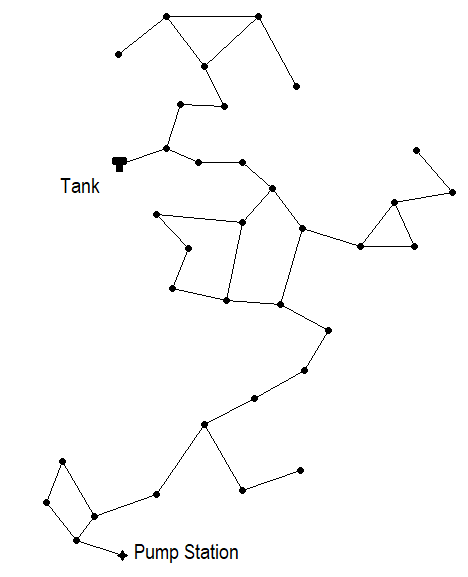}
	\caption{EPANET Example Network-2 of a WDS from Cherry Hills, CT~\cite{Rossman94net2}.}
	\label{fig:net2map}
\end{figure}

\begin{figure}[t]
	\centering
	\includegraphics[scale=0.2]{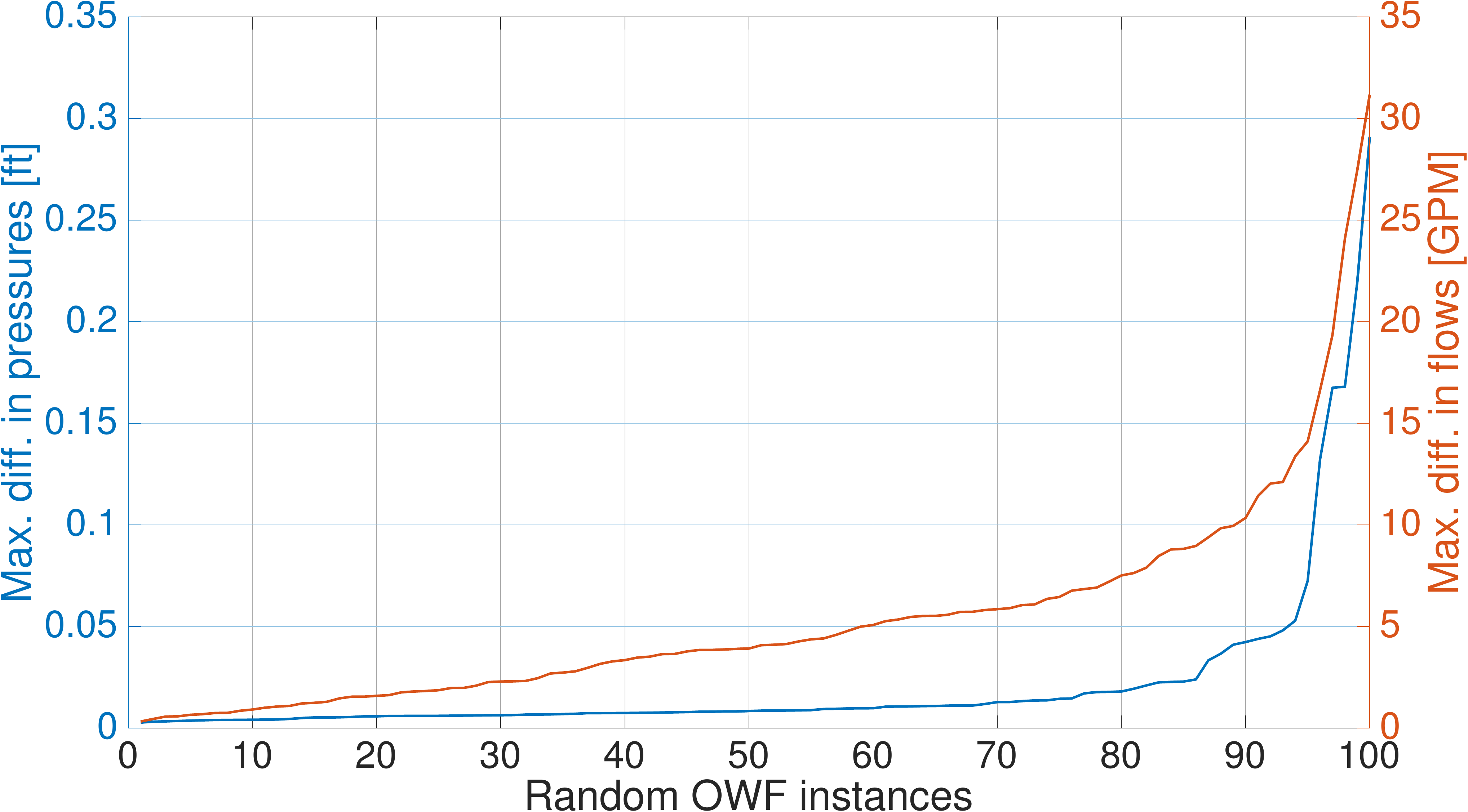}
	\caption{Maximum errors in nodal pressures and pipeline flows in the minimizer of \eqref{eq:P3} obtained for the WDS of Fig.~\ref{fig:net2map}.}
	\label{fig:hfgap}
\end{figure}

The feasibility of a minimizer obtained from \eqref{eq:P3} was also evaluated on the EPANET \emph{Example Network-2} representing a WDS from Cherry Hills, Connecticut~\cite{Rossman94net2}, which is shown in Fig.~\ref{fig:net2map}. This WDS consists of $40$ pipes, $34$ demand nodes, one tank and one pump station. Observe that none of the cycles in this WDS satisfy the assumptions of Lemma~\ref{le:parallel}. We modified the network by representing the pump station as a reservoir with pressure $100$ft connected to a fixed-speed pump with a head gain of $100$~ft. Assuming all nodes to be at the same reference elevation, the minimum pressure requirement for all nodes was set to $90$ft. The pipe friction coefficients $c_{mn}$'s, tank dimensions and the base nodal demands $d_m$'s were derived from the related EPANET file. 

To empirically evaluate the feasibility of a minimizer of \eqref{eq:P3}, we generated $100$ triplets of hourly nodal demands upon scaling the base demand by an independent uniform random variable within $[0,1]$. These hourly demands were used to solve $100$ instances of the OWF problem on a horizon of $T=3$ hours with $\lambda=10$. The maximum value of $|h_m^t-h_n^t|-c_{mn}(d_{mn}^t)^2$ for all pipes and times was recorded for all $100$ instances. These values were found to lie within $[8\cdot 10^{-5},0.56]$ with their median at $0.017$. To further understand the physical feasibility of the obtained minimizers, the nodal demands, tank injections, pump status, and reservoir pressures were used to solve a water flow (WF) problem to find the resulting nodal pressures and pipeline flows. A constrained energy function minimization-based WF solver was used from~\cite{SinghKekatosWF19}. The true pipeline flows and nodal pressures obtained from the WF solver were then compared to the corresponding values from the minimizers of \eqref{eq:P3} to quantify the error. The ranked maximum absolute differences in nodal pressures and pipeline flows for the $100$ problem instances are shown in Fig.~\ref{fig:hfgap}. Considering that the nodal pressures are around $90-190$ft and network demands are in the order of $200$~GPM, the feasibility gap for a minimizer of \eqref{eq:P3} is small for a large number of problem instances. Specifically, in $90\%$ of the instances, the maximum error in computed pressures was less than $0.04$ft, while the maximum error in computed flows was less than $10.3$~GPM.

On the computational side, the running times for the $100$ OWF instances lied in the range of $[7.5,39.1]$sec, with their median at $39$sec. The time horizon was limited to $T=3$ to reduce the running time and focus on the feasibility of \eqref{eq:P3}. Observe that MI-SOCP problems are hard in general, their computational complexity is not polynomial with respect to the number of variables and constraints, and it may change significantly across problem instances.

\begin{figure}[t]
	\centering
	\includegraphics[scale=0.28]{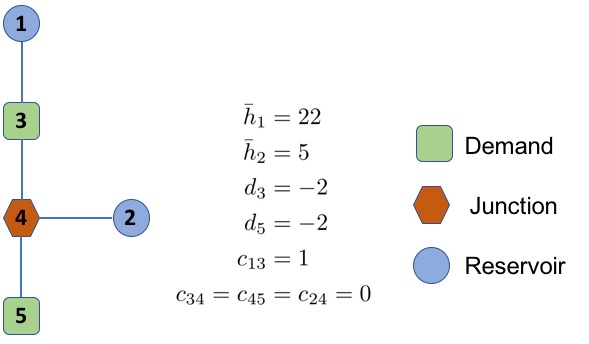}
	\caption{A simple WDS for which the relaxation is inexact.}
	\label{fig:fail}
\end{figure} 

\begin{table}[t]
	\centering
	\caption{Inexact Relaxation for the WDS of Fig.~\ref{fig:fail}}
	\vspace*{-1em}
	\begin{tabular}{|c| c| c|c|} 
		\hline\hline
		Variable&\eqref{eq:P3}&OWF in \cite{TaylorOWFCDC}&\eqref{eq:P1}\\
		\hline
		$h_1$&10&10&22\\\hline
		$h_2$&5&5&6\\\hline
		$h_3$&6&6&6\\\hline
		$h_4$&5&5&6\\\hline
		$h_5$&5&5&6\\\hline
		$d_{13}$&2&2&4\\\hline
		$d_{34}$&0&0&2\\\hline
		$d_{24}$&2&2&0\\\hline
		$d_{45}$&2&2&2\\\hline
		\emph{comment}& inexact & inexact & optimal\\\hline
	\end{tabular}
	\label{tbl:inexact}
\end{table}

Finally, to provide an example of inexact relaxation, we built the WDS of Figure~\ref{fig:fail}. Problem~\eqref{eq:P3} and the OWF scheme of \cite{TaylorOWFCDC} were solved on this WDS for minimum pressures at nodes $3$, $4$, and $5$, set to $6$, $0$, and $0$. This setup features a unique feasible point: Since all edges but $(1,3)$ are lossless, nodes $2-5$ must have equal pressures. Because $\underline{h}_3=6$m, the second reservoir with $\bar{h}_2=5$m cannot supply water, the entire demand must be fulfilled by reservoir $1$. This feasible point is shown in Table~\ref{tbl:inexact}, along with the minimizers of \eqref{eq:P3} and \cite{TaylorOWFCDC}. Both relaxed schemes yielded an infeasible point for \eqref{eq:P1}. The solver of \cite{TaylorOWFCDC} was not tested on the $10$-node WDS earlier because it presumes: \emph{i)} variable-speed pumps with speeds that can reach zero; and \emph{ii)} that once a solution $(\tbd,\bd)$ is found, a feasible pressure $\bh$ can always be obtained. 

\section{Conclusions}\label{sec:conclusions}
To cater a more adaptive WDS operation, optimal pump scheduling has been formulated here as an OWF task. Different from existing formulations, the developed OWF model includes critical pressure constraints capturing the operation of tanks, reservoirs, pipes, and valves. The original mixed-integer non-convex problem has been modified to a mixed-integer second-order cone program over a relaxed feasible set. Moreover, its objective is augmented by a judiciously designed penalty term, so that under specific conditions, this modified problem formulated as an MI-SOCP can recover minimizers of the original problem. Numerical tests validate that by properly tuning the penalization parameter $\lambda$, the modified problem solves the original OWF over different scenarios of water demand and electricity pricing. 

Off-the-shelf MI-SOCP solvers have improved significantly over the last years, yet MI-SOCP's bear no computational complexity guarantees. Although a related MI-SOCP-based solver we have developed in~\cite{SinghKekatosWF19} for the water flow problem scales well with the network size, that is not always the case here for \eqref{eq:P3}. The running time of \eqref{eq:P3} depends on water demands, electricity prices, and the values of $M$'s involved in the big-$M$ constraints. To accelerate \eqref{eq:P3}, future research could pursue two directions. First, one could exploit the temporal dynamics of OWF. Water system decisions are coupled across time only through the tank operation of \eqref{eq:tank1}. Therefore, one could select tank levels $\{\ell_m^t\}_{m\in\mcM_b}$ as the system states; discretize their values based on the desired approximation/complexity trade-off; and handle \eqref{eq:P3} using approximate dynamic programing. Secondly, based on prior experience, the WDS operator may be able to fix some of the binary variables capturing the flow directions on pipes and the operating statuses of pumps/reservoirs, to prespecified values. 

Other pertinent research directions include generalizing our OWF formulation towards scheduling variable-speed pumps and/or incorporating stochasticity in water demands and electricity prices. Finally, the developed framework could be readily used for jointly scheduling WDS and electric power distribution networks to realize the vision for smart cities.

\appendix
\begin{IEEEproof}[Proof of Theorem~\ref{th:1}]
Being a minimizer, $\tbs_3$ is also feasible for \eqref{eq:P3}. A feasible point of \eqref{eq:P3} satisfies only those equations in \eqref{eq:heq} related to pumps. The equality constraints in \eqref{eq:heq} corresponding to lossy pipes are replaced by one-sided linear \emph{inequality} constraints in \eqref{eq:P3}. To express these facts in a matrix-vector notation, partition $\bA(\bd^t)$ into submatrix $\bA_p(\bd^t)$ having the rows of $\bA(\bd^t)$ related to pumps; and submatrix $\bA_l(\bd^t)$ having the rows related to lossy pipes. The rows of $\bA(\bd^t)$ can be permuted without loss of generality so that
	\begin{equation}\label{eq:A}
	\bA(\bd^t)=\left[\begin{array}{c}
	\bA_p(\bd^t)\\
	\bA_l(\bd^t)
	\end{array}\right].
	\end{equation}
Likewise, the mapping $\bb(\tbd^t,\bd^t)$ in \eqref{eq:heq} can be partitioned into $\bb_p(\tbd^t)$ and $\bb_l(\bd^t)$. A vector $\bh$ is feasible for the relaxed problem \eqref{eq:P3} if instead of \eqref{eq:heq}, it satisfies
\begin{subequations}\label{eq:th1LP}
	\begin{align}
	\bA_p(\bd^t)\bh^t&=\bb_p(\tbd^t)\label{seq:th1LPa},\quad\forall t\\
	\bA_l(\bd^t)\bh^t&\geq\bb_l(\bd^t)\label{seq:th1LPb}\geq \mathbf{0},\quad\forall t.	
	\end{align}
\end{subequations}

Granted $\mcH(\tbd_3,\bd_3)$ is non-empty by hypothesis, there exists an $\cbh_3\in\mcH(\tbd_3,\bd_3)$ so that $\cbs_3:=\{\tbd_3,\bd_3,\cbh_3\}$ satisfies \eqref{eq:heq}--\eqref{eq:hineq}. Because $\cbs_3$ satisfies \eqref{eq:heq}, it satisfies the constraints \eqref{seq:th1LPb} with equality. Thus, vector $\cbs_3$ is feasible for \eqref{eq:P3}. Moreover, the cost of \eqref{eq:P3} for $\cbs_3$ is $f(\tbd_3)+\lambda g(\cbh_3)=f_3+\lambda\sum_{t=1}^{T} \|\bA_l(\bd_3^t)\cbh_3^t\|_1$, where $f_3:=f(\tbd_3)$, and $\cbh_3^t$ and $\bd_3^t$ are accordingly the subvectors of $\cbh_3$ and $\bd_3$ collecting the entries corresponding to time $t$. Since $\cbs_3$ satisfies \eqref{seq:th1LPb} with equality, the cost becomes $f_3+\lambda\sum_{t=1}^{T} \|\bb_l(\bd_3^t)\|_1$.
	
Proving by contradiction, suppose $\bh_3\notin\mcH(\tbd_3,\bd_3)$. This implies $\bh_3$ does not satisfy the left-hand side of \eqref{seq:th1LPb} with equality. Instead, there exists a sequence of $\bepsilon^t\geq \mathbf{0}$, such that $\bA_l(\bd_3^t)\bh_3^t=\bb_l(\bd_3^t)+\boldsymbol{\epsilon}^t$ for all $t$ and $\sum_{t=1}^{T}\boldsymbol{\epsilon}^t\neq \mathbf{0}$. Evaluating the objective of \eqref{eq:P3} for the minimizer $\bs_3$ yields
\begin{align*}
	f(\tbd_3)+\lambda g(\tbh_3)&=f_3+\lambda\sum_{t=1}^{T}\|\bA_l(\bd_3^{t})\bh_3^{t}\|\nonumber\\
	&=f_3+\lambda\sum_{t=1}^{T}\left(\|\bb_l^t(\bd_3^t)\|_1+\|\boldsymbol{\epsilon}_t\|_1\right)\nonumber\\
	&>f_3+\lambda\sum_{t=1}^{T}\|\bb_l^t(\bd_3^t)\|_1
\end{align*} 
where the second equality stems from $\bb_l(\bd^t)\geq \mathbf{0}$ and $\boldsymbol{\epsilon}_t\geq \mathbf{0}$ for all $t$; and the strict inequality holds because $\lambda>0$ and $\sum_{t=1}^{T}\boldsymbol{\epsilon}^t\neq \mathbf{0}$. This inequality contradicts the optimality of $\bs_3$, and nullifies the hypothesis that $\bh_3\notin\mcH(\tbd_3,\bd_3)$.
\end{IEEEproof}

\begin{IEEEproof}[Proof of Lemma~\ref{le:parallel}]
Since this proof refers to a particular time, the superscript $t$ is omitted for simplicity. Given a point $\{\tbd,\bd,\bh\}$, an edge will be termed (in)exact if constraint \eqref{eq:relaxed} is satisfied with (in)equality for that point. Since all nodes incident to $\mcR$ excluding $m$ host no tanks or reservoirs, they must have non-positive injections. Therefore, its two incident edges cannot both have outgoing water flows from \eqref{eq:node1}. This implies that the ring can either consist of two parallel paths, or a directed cycle. In the latter case, adding the constraints $h_{i}-h_{j}\geq c_{ij}(d_{ij})^2$ around $\mcR$ would give $\sum_{(i,j)\in \mcR}c_{ij}d_{ij}^2\leq h_m-h_m=0$, implying $d_{ij}=0$ for all edges in $\mcR$, which is a contradiction. Thus, the ring $\mcR$ consists of two parallel paths from $m$ to some node $n$, henceforth termed $\mcP_1$ and $\mcP_2$. 

The rest of the proof proceeds in two steps. The first step shows there exists a minimizer of \eqref{eq:P3} with at most one inexact edge in $\mcR$. The second step reduces the number to none. 

For the first step, we will modify the pressure vector in $\bs_3$ to construct $\hbs_3:=\{\tbd_3,\bd_3,\hat{\bh}_3\}$ for which there exists at most one inexact edge in $\mcR$. The new point $\hbs_3$ is feasible for \eqref{eq:P3} and attains smaller or equal cost than $\bs_3$. To do so, for each node $k$ incident to $\mcR$ excluding $m$ and $n$, assign the pressure consistent with \eqref{eq:headloss} along the path $\mcP_{mk}$ from $m$ to $k$:
\[\hat{h}_k:=h_m-\sum_{(i,j)\in \mcP_{mk}}c_{ij}d_{ij}^2\geq h_k\geq \underline{h}\]
where the first inequality stems from summing up the constraints $h_{i}-h_{j}\geq c_{ij}d_{ij}^2$ for all edges $(i,j)$ along $\mcP_{mk}$, and guarantees that $\hat{h}_k$ is feasible. 

For the terminal node $n$, assign the pressure
\begin{equation}\label{eq:hhn:1}
\hat{h}_n:=\min_{l\in\{1,2\}}\Big\{h_m-\sum_{(i,j)\in \mcP_l}c_{ij}d_{ij}^2\Big\}.
\end{equation}
Adding the constraints $h_{i}-h_{j}\geq c_{ij}d_{ij}^2$ for all edges $(i,j)$ in $\mcP_l$ and $\mcP_2$ separately, yields
\begin{equation}\label{eq:hhn1}
 h_m-h_n\geq \sum_{(i,j)\in \mcP_l}c_{ij}d_{ij}^2,\quad l\in\{1,2\}.
\end{equation}
Hence, we get that
\begin{equation}\label{eq:hhn2}
h_n\leq \min_{l\in\{1,2\}}\Big\{h_m-\sum_{(i,j)\in \mcP_l}c_{ij}d_{ij}^2\Big\}=\hat{h}_n
\end{equation}
implying $\hat{h}_n\geq h_n\geq \underline{h}$. 

Since the pressures on the nodes within $\mcR$ have been increased and they are not upper bounded in the absence of tanks or reservoirs, the point $\hbs_3$ is feasible. The difference in the objective of \eqref{eq:P3} attained by $\bs_3$ and $\hbs_3$ is
\begin{align*}
&f(\tbd_3)+\lambda g(\bh_3)-f(\tbd_3)-\lambda g(\hbh_3)\\
&=\lambda\sum_{(i,j)\in \mcR}\left(|h_i-h_j|-|\hat{h}_i-\hat{h}_j|\right).
\end{align*}
Since all directed edges in $\mathcal{P}_1$ and $\mathcal{P}_2$ have positive flows
\begin{align*}
\sum_{(i,j)\in \mcR}|h_i-h_j|&=\sum_{(i,j)\in \mcP_1}(h_i-h_j)+\sum_{(i,j)\in \mcP_2}(h_i-h_j)\nonumber\\
&=2(h_m-h_n).
\end{align*}
Applying the same argument for $\hat{\bh}_3$, it follows that
\begin{equation}\label{eq:hhn3}
f(\tbd_3)+\lambda g(\bh_3)-f(\tbd_3)-\lambda g(\hbh_3)=2\lambda (\hat{h}_n-h_n)\geq 0.
\end{equation}

If for $\bs_3$ there exist inexact edges in both $\mcP_1$ and $\mcP_2$, then \eqref{eq:hhn1} holds with strict inequality for both paths. It follows from \eqref{eq:hhn2} that $\hat{h}_n>h_n$, and so $\hbs_3$ contradicts the optimality of $\bs_3$. This proves that all inexact edges in $\mcR$ must belong exclusively to $\mcP_1$ or $\mcP_2$. In the latter case, the inequality in \eqref{eq:hhn2} holds with equality, and from \eqref{eq:hhn3} the point $\hbs_3$ becomes a minimizer of \eqref{eq:P3}. Note $\hbs_3$  has at most one inexact edge in $\mcR$, and that is the last edge in $\mcP_1$ or $\mcP_2$. 

For the second step of this proof and proving by contradiction, suppose there exist exactly one inexact edge for the minimizer $\bs_3$ in $\mcP_1$. That means that \eqref{eq:hhn1} holds with inequality for $l=1$, and equality for $l=2$, implying
\begin{align}\label{eq:le:1}
\sum_{(i,j)\in \mcP_1}c_{ij}d_{ij}^2<\sum_{(i,j)\in \mcP_2}c_{ij}d_{ij}^2.
\end{align}

From $\bd_3$, construct a water flow vector $\check{\bd}_3$ with entries
\begin{equation}\label{eq:redist}
\check{d}_{ij}=
\left\{\begin{array}{ll}
	d_{ij}+\epsilon &,~(i,j)\in\mcP_1\\
	d_{ij}-\epsilon &,~(i,j)\in\mcP_2\\
	d_{ij} &,~(i,j)\in\mcP\setminus(\mcP_1\cup\mcP_2)
\end{array}\right.
\end{equation}  
for some $\epsilon>0$. This redistribution of flows satisfies \eqref{eq:node1}. Moreover, for increasing $\epsilon$, the LHS of \eqref{eq:le:1} increases and the RHS decreases. This is because $c_{ij}d_{ij}^2$ is an increasing function for positive $d_{ij}$. The goal is to select $\epsilon$, so that 
\begin{equation}\label{eq:le1:equality}
\sum_{(i,j)\in \mcP_1}c_{ij}\check{d}_{ij}^2=\sum_{(i,j)\in \mcP_2}c_{ij}\check{d}_{ij}^2<\sum_{(i,j)\in \mcP_2}c_{ij}d_{ij}^2.
\end{equation}
While increasing $\epsilon$ to achieve \eqref{eq:le1:equality}, some of the $\{\check{d}_{ij}\}_{(i,j)\in\mcP_2}$ may become negative. This case is ignored for now. 

Construct next a new pressure vector $\cbh_3$ by changing the entries of $\bh_3$ corresponding to the non-root nodes in $\mcR$ as
\begin{align}\label{eq:le1:ch}
\check{h}_k:=h_m-\sum_{(i,j)\in \mcP_{mk}}c_{ij}\check{d}_{ij}^2.
\end{align}
For $k=n$, the sum in the RHS of \eqref{eq:le1:ch} can be evaluated over $\mcP_1$ or $\mcP_2$, since these two sums are equal from \eqref{eq:le1:equality}. The constructed pressures for nodes incident to $\mcR$ satisfy
\begin{equation}\label{eq:le1:sequence}
\check{h}_k\geq \check{h}_n> h_n \geq \underline{h}.
\end{equation}
The first inequality holds because node $n$ has the largest value for the sum in \eqref{eq:le1:ch}; and the second inequality because
\[\check{h}_n=h_m-\sum_{(i,j)\in \mcP_2}c_{ij}\check{d}_{ij}^2>h_m-\sum_{(i,j)\in \mcP_2}c_{ij}d_{ij}^2=h_n.\]
The inequalities in \eqref{eq:le1:sequence} prove that $\cbh_3$, and hence the point $\cbs_3:=\{\tbd_3,\cbd_3,\cbh_3\}$ is feasible for \eqref{eq:P3}. The difference in the objective of \eqref{eq:P3} attained by $\bs_3$ and $\cbs_3$ is
\begin{align*}
f(\tbd_3)+\lambda g(\bh_3)-f(\tbd_3)-\lambda g(\cbh_3)=2\lambda (\check{h}_n-h_n)>0
\end{align*}
which contradicts the optimality of $\bs_3$.

Since all water injections at non-root nodes over $\mcR$ are non-positive, the water flows are non-increasing along $\mcP_2$. This implies that $d_{ij}\geq d_{n_1,n}$ for all $(i,j)\in\mcP_2$, where $(n_1,n)$ is the last edge of $\mcP_2$. Thus, by increasing $\epsilon$, the flow $d_{n_1,n}$ may become negative. In that case, the edge $(n_1,n)$ is removed from $\mcP_2$ and appended to $\mcP_1$, forming a new pair of parallel paths with $n_1$ as the new terminal node. The second step of this proof can be repeated on the new parallel paths.
\end{IEEEproof}

\begin{IEEEproof}[Proof of Theorem~\ref{th:2}]
Let $\mcT:= (\mcM,\mcP_\mcT)$ be a spanning tree of $(\mcM,\mcP(\bd_3^t))$. Reorder the equations in \eqref{eq:heq} as
\begin{equation}\label{eq:partition}
\left[\begin{array}{c}
	\bA_\mcT(\bd_3^t)\\
	\bA_{\bar{\mcT}}(\bd_3^t)
	\end{array}\right]\bh^t=\left[\begin{array}{c}
	\bb_\mcT(\tbd_3^t,\bd_3^t)\\
	\bb_{\bar{\mcT}}(\tbd_3^t,\bd_3^t)
	\end{array}\right]
\end{equation}
where $\bA_\mcT(\bd_3^t)$ and $\bb_\mcT(\tbd_3^t,\bd_3^t)$ are the rows of $\bA(\bd_3^t)$ and $\bb(\tbd_3^t,\bd_3^t)$ corresponding to the edges in $\mcP_\mcT$; and  $\bA_{\bar{\mcT}}(\bd_3^t)$ and $\bb_{\bar{\mcT}}(\tbd_3^t,\bd_3^t)$ the rows corresponding to the edges in $\mcP\setminus\mcP_\mcT$. 

Being an edge-node incidence matrix for a tree, matrix $\bA_\mcT(\bd_3^t)$ is full row-rank~\cite{GodsilRoyle}, and hence the system $\bA_\mcT(\bd_3^t)\bh^t=\bb_\mcT(\tbd_3^t,\bd_3^t)$ is consistent. The rows of $\bA_{\bar{\mcT}}(\bd_3^t)$ correspond to the links defined by $\mcT$. By the hypothesis, every undirected cycle in $\left(\mcM,\mcP(\bd_3^t)\right)$ is a ring. Then, all but one of its edges belong to $\mcT$, and the remaining edge belongs to $\bar{\mcT}$. In fact, every edge in $\bar{\mcT}$ must belong to a ring. Since by the conditions of Lemma~\ref{le:parallel}, no pumps are allowed on a ring, every equation in the bottom part of \eqref{eq:partition} corresponds to a lossy pipeline $(k,l)$ and will be of the form $h_k^t-h_l^t=c_{kl}d_{kl}^2$.

Since we refer to time $t$, the superscript $t$ is omitted to unclutter notation. Consider link $(k,l)\in\bar{\mcT}$ that belongs to the pair of parallel paths $\mcP_1$ and $\mcP_2$ with origin node $m$ and destination $n$. Without loss of generality, let also $(k,l)\in\mcP_1$. From Lemma \ref{le:parallel}, it holds that $h_i-h_j=c_{ij}d_{ij}^2$ for all $(i,j)\in\mcP_1\cup\mcP_2$. Summing these constraints along $\mcP_1$ and $\mcP_2$ yields
\begin{subequations}\label{eq:equal}
\begin{align}
\sum_{(i,j)\in \mcP_1}(h_i-h_j)&=\sum_{(i,j)\in \mcP_1}c_{ij}d_{ij}^2=h_m-h_n\label{eq:equal:P1}\\
\sum_{(i,j)\in \mcP_2}(h_i-h_j)&=\sum_{(i,j)\in \mcP_2}c_{ij}d_{ij}^2=h_m-h_n\label{eq:equal:P2}
\end{align}
\end{subequations}
so that \eqref{eq:equal:P1} equals \eqref{eq:equal:P2}. Separating the contribution of edge $(k,l)$ from $\mcP_1$ in the leftmost and central parts of \eqref{eq:equal:P1} provides
\begin{subequations}\label{eq:kl}
\begin{align}
 h_k-h_l&=\sum_{(i,j)\in \mcP_2}(h_i-h_j)-\sum_{(i,j)\in \mcP_1\setminus (k,l)}(h_i-h_j)\\
 c_{kl}d_{kl}^2&=\sum_{(i,j)\in \mcP_2}c_{ij}d_{ij}^2-\sum_{(i,j)\in \mcP_1\setminus(k,l)}c_{ij}d_{ij}^2.
\end{align}
\end{subequations}
Note that the pressure drop equations along for all edges $(i,j)\in\mcP_1\cup\mcP_2\setminus(k,l)$ are rows in the system $\bA_\mcT(\bd_3^t)\bh^t=\bb_\mcT(\tbd_3^t,\bd_3^t)$. From \eqref{eq:kl}, the pressure drop equation corresponding to edge $(k,l)\in\bar{\mcT}$ has been expressed as a linear combination of the rows of $\bA_{\mcT}(\bd_3)\bh= \bb_{\mcT}(\tbd_3,\bd_3)$. The argument holds for all equations in the bottom part of \eqref{eq:partition}, thus making the overall system in \eqref{eq:heq} consistent.
\end{IEEEproof}

\balance
\bibliography{myabrv,water}
\bibliographystyle{IEEEtran}

\begin{IEEEbiography}[{\includegraphics[width=1in,height=1.25in,clip,keepaspectratio]{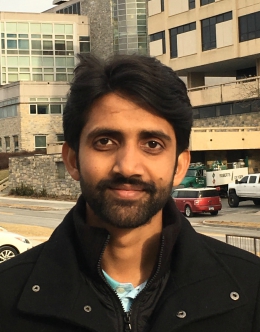}}] {Manish K. Singh} received the B.Tech. degree from the Indian Institute of Technology (BHU), Varanasi, India, in 2013;  and the M.S. degree from Virginia Tech, Blacksburg, VA, USA, in 2018; both in electrical engineering. During 2013-2016, he worked as an Engineer in the Smart Grid Dept. of POWERGRID, the central transmission utility of India. He is currently pursuing a Ph.D. degree at Virginia Tech. His research interests are focused on the application of optimization, control, and graph-theoretic techniques to develop algorithmic solutions for operation and analysis of water, natural gas, and electric power systems.
\end{IEEEbiography} 

\begin{IEEEbiography}[{\includegraphics[width=1in,height=1.25in,clip,keepaspectratio]{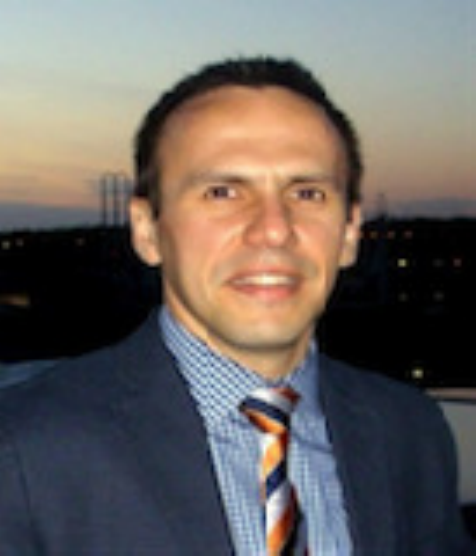}}] {Vassilis Kekatos} (SM'16) is an Assistant Professor with the Bradley Dept. of ECE at Virginia Tech. He obtained his Diploma, M.Sc., and Ph.D. from the Univ. of Patras, Greece, in 2001, 2003, and 2007, respectively. He is a recipient of the NSF Career Award in 2018 and the Marie Curie Fellowship. He has been a research associate with the ECE Dept. at the Univ. of Minnesota, where he received the postdoctoral career development award (honorable mention). During 2014, he stayed with the Univ. of Texas at Austin and the Ohio State Univ. as a visiting researcher. His research focus is on optimization and learning for future energy systems. He is currently serving in the editorial board of the IEEE Trans. on Smart Grid.
\end{IEEEbiography}

\end{document}

%% file: stylesV2.tex
\DeclareMathOperator{\sign}{sign}

\newtheorem{lemma}{Lemma}
\newtheorem{corollary}{Corollary}
\newtheorem{theorem}{Theorem}

\newcommand \bb{\mathbf{b}}

\newcommand \bd{\mathbf{d}}

\newcommand \bh{\mathbf{h}}

\newcommand \bs{\mathbf{s}}

\newcommand \bx{\mathbf{x}}

\newcommand \bA{\mathbf{A}}

\newcommand \bepsilon{\boldsymbol{\epsilon}}

\newcommand \mcG{\mathcal{G}}
\newcommand \mcH{\mathcal{H}}

\newcommand \mcM{\mathcal{M}}

\newcommand \mcP{\mathcal{P}}

\newcommand \mcR{\mathcal{R}}
\newcommand \mcS{\mathcal{S}}
\newcommand \mcT{\mathcal{T}}

\newcommand \mcX{\mathcal{X}}

\newcommand \bmcP{\bar{\mathcal{P}}}

\newcommand \tbd{\tilde{\mathbf{d}}}

\newcommand \tbh{\tilde{\mathbf{h}}}

\newcommand \tbs{\tilde{\mathbf{s}}}

\newcommand \cbd{\check{\mathbf{d}}}

\newcommand \cbh{\check{\mathbf{h}}}

\newcommand \cbs{\check{\mathbf{s}}}

\newcommand \hbh{\hat{\mathbf{h}}}

\newcommand \hbs{\hat{\mathbf{s}}}
